\magnification=\magstep1
\input amstex
\documentstyle{amsppt}

\define\defeq{\overset{\text{def}}\to=}

\define\ab{\operatorname{ab}}
\define\Aut{\operatorname{Aut}}
\define\cl{\operatorname{cl}}
\define\et{\operatorname{et}}
\define\fin{\operatorname{fin}}
\define\Gal{\operatorname{Gal}}
\define\OSub{\operatorname{OSub}}
\define\pr{\operatorname{pr}}
\define\rank{\operatorname{rank}}
\define\sep{\operatorname{sep}}
\define\Sub{\operatorname{Sub}}

\define\tr{\operatorname{tr}}
\define\ur{\operatorname{ur}}
\def \ord{\operatorname {ord}}
\def \div{\operatorname {div}}
\def \Div{\operatorname {Div}}
\def \Pri{\operatorname {Pri}}

\def \c{\operatorname {c}}

\def \ab{\operatorname {ab}}
\def \isom {\overset \sim \to \rightarrow}
\def \Ker{\operatorname {Ker}}
\def \Hom{\operatorname {Hom}}
\def \char{\operatorname {char}}

\define\Primes{\frak{Primes}}

\NoRunningHeads
\NoBlackBoxes
\topmatter

\title
A Prime-to-$p$ Version of Grothendieck's Anabelian Conjecture
for Hyperbolic Curves over Finite Fields of Characteristic $p>0$
\endtitle


\author
Mohamed Sa\"\i di and Akio Tamagawa
\endauthor





\abstract
In this paper, we prove a prime-to-$p$ version of Grothendieck's anabelian
conjecture for hyperbolic curves over finite fields of characteristic $p>0$, 
whose original (full profinite) version was proved by Tamagawa in the affine 
case and by Mochizuki in the proper case.
\endabstract

\toc

\subhead
\S 0. Introduction
\endsubhead

\subhead
\S 1. Characterization of Decomposition Groups
\endsubhead

\subhead
\S 2. 
Cuspidalizations of Proper Hyperbolic Curves
\endsubhead

\subhead
\S 3. Kummer Theory and Anabelian Geometry
\endsubhead

\subhead	
\S 4. Recovering the Additive Structure
\endsubhead

\endtoc

\endtopmatter

\document            

\subhead
\S 0. Introduction
\endsubhead

Let $k$ be a finite field of characteristic $p>0$ and 
$U$ a hyperbolic curve over $k$. Namely, $U=X-S$, where 
$X$ is a proper, smooth, geometrically connected curve 
of genus $g$ over $k$ and $S\subset X$ is a divisor 
which is finite \'etale of degree $r$ over $k$, such that 
$2-2g-r<0$. 
We have the following exact sequence of profinite groups:
$$1\to \pi_1(U\times _k{\bar k},*)\to \pi_1(U,*)\to G_k\to 1.$$
Here, $G_k$ is the absolute Galois group $\Gal (\bar k/k)$, $*$ 
means a suitable geometric 
point, and $\pi _1$ stands for the \'etale fundamental group. 
The following result is fundamental
in the anabelian geometry of hyperbolic curves over finite fields.

\proclaim{Theorem (Tamagawa, Mochizuki)} Let $U$, $V$ be 
hyperbolic curves over 
finite fields $k_U$, $k_V$, respectively. Let 
$$\alpha : \pi_1(U,*)\isom \pi_1(V,*)$$
be an isomorphism of profinite groups. 
Then $\alpha$ arises from a uniquely determined 
commutative diagram of schemes:
$$
\CD
\Tilde U@>{\sim}>> \Tilde V\\
@VVV   @VVV \\
U @>{\sim}>> V\\
\endCD
$$
in which the horizontal arrows are isomorphisms, and the vertical arrows are the 
profinite \'etale universal coverings determined 
by the profinite groups $\pi_1(U,*)$, $\pi_1(V,*)$, respectively.
\endproclaim

This result was proved by Tamagawa (cf. [Tamagawa1], Theorem (4.3)) in the affine case 
(together with a certain tame version), 
and more recently by Mochizuki (cf. [Mochizuki2], Theorem 3.12) in the proper case. 
It implies, in particular, that one can embed a suitable category of hyperbolic 
curves over finite fields into the category of profinite groups. It is essential in 
the anabelian philosophy of Grothendieck, as was formulated in [Grothendieck], 
to be able to determine the image of this functor. Recall that the full structure 
of the profinite group $\pi_1(U\times _k{\bar k},*)$ is unknown 
(for any single example of 
$U$ which is hyperbolic). Hence, a fortiori, the structure of  
$\pi_1(U,*)$ is unknown. 
(Even if we replace the fundamental groups 
$\pi_1(U\times _k{\bar k},*)$, $\pi_1(U,*)$ by the tame fundamental groups 
$\pi_1^t(U\times _k{\bar k},*)$, $\pi_1^t(U,\ast)$, respectively, 
the situation is just the same.) 
Thus, the problem of determining the image of the above functor 
seems to be quite difficult, at least for the moment. 
In this paper we investigate the following question:

\definition{Question 1} 
Is it possible to prove 
any result analogous to the above Theorem 
where $\pi_1(U,*)$ is replaced by some (continuous) quotient of 
$\pi_1(U,*)$ whose structure is better understood?
\enddefinition

The first quotients that come into mind are the following. Let $\Cal C$ (respectively, 
$\Cal C^l$) be the class of finite groups of order prime to $p$ (respectively, 
finite $l$-groups, where $l\neq p$ is a fixed prime number). 
Let $\Delta _U$ be the maximal pro-prime-to-$p$ 
(i.e., pro-$\Cal C$) quotient of $\pi_1(U\times _k{\bar k},*)$. 
For a profinite group $\Gamma$, $\Gamma^l$ stands for the 
maximal pro-$l$ (i.e., pro-$\Cal C^l$) quotient of $\Gamma$. 
Thus, in particular, 
$\Delta _U^l$ coincides with 
$\pi_1(U\times _k{\bar k},*)^l$. 
Here, the structures of $\Delta_U$ and $\Delta_U^l$ are well understood 
--- $\Delta_U$ (respectively, $\Delta_U^l$) is isomorphic to 
the pro-prime-to-$p$ (respectively, pro-$l$) completion of a 
certain well-known finitely generated discrete group (i.e., either a free 
group or a surface group). 
Let $\Pi_U\overset \text {def}\to=\pi_1(U,*)/\Ker(\pi_1
(U\times _k{\bar k},*)\twoheadrightarrow \Delta_U)$ (respectively, $\Pi_U^{(l)}
\defeq \pi_1(U,*)/\Ker 
(\pi_1(U\times _k{\bar k},*)\twoheadrightarrow \Delta _U^l)$ 
be the 
corresponding 
quotients of $\pi_1(U\times _k{\bar k},*)$. We shall refer to $\Pi_U$ 
as the geometrically pro-$\Sigma_U$ \'etale fundamental group of $U$, where 
$\Sigma_U\defeq\Primes - 
\{\char (k)\}$, and $\Primes$ stands for the set of all prime numbers.

\definition{Question 2} 
Is it possible to prove any result analogous to the above Theorem 
where $\pi_1(U,*)$ is replaced by $\Pi _U$, $\Pi _U^{(l)}$, respectively?
\enddefinition

Our main result in this paper is the following 
%
%
%
%
(cf. Corollary 3.10):

\proclaim{Theorem 1 (A Prime-to-$p$ Version of Grothendieck's 
Anabelian Conjecture for 
Hyperbolic Curves over Finite Fields)} 
Let $U$, $V$ be 
hyperbolic curves over finite fields $k_U$, $k_V$, 
respectively. 
Let $\Sigma _U \overset \text {def}\to=\Primes-\{\char (k_U)\}$,  
$\Sigma _V\overset \text {def}\to=\Primes-\{\char (k_V)\}$,
and write $\Pi_U$, $\Pi_V$ for the geometrically pro-$\Sigma _U$ 
\'etale 
fundamental group of $U$, and the geometrically pro-$\Sigma _V$ 
\'etale 
fundamental group of $V$, respectively.
Let
$$\alpha:\Pi_U\isom \Pi_V$$
be an isomorphism of profinite groups. Then $\alpha$ arises from a uniquely 
determined commutative diagram of schemes:
$$
\CD
\Tilde U @>{\sim}>> \Tilde V \\
@VVV   @VVV \\
U @>{\sim}>> V \\
\endCD
$$
in which the horizontal arrows are isomorphisms and the vertical arrows are the 
profinite \'etale coverings corresponding to the groups $\Pi_U$, $\Pi_V$, 
respectively. 
\endproclaim

As an important consequence of Theorem 1, 
we deduce in Corollary 3.11 
the following prime-to-$p$ version of Uchida's Theorem on isomorphisms 
between absolute Galois groups of function fields (cf. [Uchida]).

\proclaim{Theorem 2} Let $X$, $Y$ be proper, smooth, geometrically connected 
curves over finite fields $k_X$, $k_Y$, respectively. 
Let $K_X$, $K_Y$ be the function fields of $X$, $Y$, respectively. Let $G_{K_{X}}$,
$G_{K_{Y}}$ be the absolute Galois groups of $K_X$, $K_Y$, respectively, 
and let $G_{K_{X}}'$, 
$G_{K_{Y}}'$ be their geometrically pro-prime-to-characteristic quotients 
(cf. discussion before Corollary 3.11). Let 
$$\alpha : G_{K_{X}}'\isom G_{K_{Y}}'$$
be an isomorphism of profinite groups. Then $\alpha$ arises from a uniquely 
determined commutative diagram of field extensions:
$$
\CD
(K_X)\sptilde @>{\sim}>> (K_Y)\sptilde \\
@AAA              @AAA \\
K_{X} @>{\sim}>> K_{Y} \\
\endCD
$$
in which the horizontal arrows are isomorphisms and the vertical arrows are the 
extensions corresponding to the Galois groups $G_{K_{X}}'$,  $G_{K_{Y}}'$, respectively. 
\endproclaim

Our proof of Theorem 1 relies substantially on the methods 
of Tamagawa and Mochizuki. 
We shall explain this briefly in the case where 
$U$ is proper (hence, $U=X$). (cf. Theorem 3.9. 
The general case can be reduced to this special case.) 
Starting from $\Pi _X$,  Tamagawa's method characterizes group-theoretically
the decomposition groups at points of $X$
in $\Pi _X$. The problem of recovering the points of $X$ from the corresponding decomposition 
groups is related to the question of whether the natural map
$$X^{\cl}\to \Sub (\Pi_X)_{\Pi_{X}} 
$$from the set of closed points of $X$ 
to the set of conjugacy classes of closed subgroups of $\Pi_X$,
which maps a point $x$ to the conjugacy class of its decomposition group 
$D_x$,
is injective. This map is known to be injective in the 
full profinite case, i.e., when one 
starts from $\pi_1(X,*)$ instead of $\Pi_X$. In our case 
we are only able to prove that the 
above map is almost injective, 
i.e., injective outside a finite set $E_X\subset X^{\cl}$. 
Thus, we can only recover from $\Pi _X$ the set of points in a nonempty open subset 
$X-E_X$. 

In [Mochizuki2], Mochizuki developed 
a theory of cuspidalizations of \'etale fundamental 
groups of proper hyperbolic curves. One of the consequences 
of the main results of this theory 
is that, starting from $\Pi_X$, one can recover in a functorial way
the Kummer theory of open affine subsets 
$U_S\defeq X-S$, where $S$ is a finite set of closed 
points contained in  $X-E_X$. Using Kummer 
theory, one then recovers (up to Frobenius twist) the multiplicative group 
$\Cal O_{E_X}^{\times}$ of rational functions on $X$ whose 
divisor has a support disjoint from $E_X$. Thus, starting from an isomorphism
$$\Pi_X\isom \Pi_Y$$ 
as in Theorem 1 we can recover, up to Frobenius twist, an injective embedding between 
multiplicative groups
$$\Cal O_{E_X}^{\times} \hookrightarrow \Cal O_{E_Y}^{\times}.$$ 
The issue is then to show that this embedding arises from 
a uniquely determined embedding of fields $K_X\hookrightarrow K_Y$, between the 
corresponding function fields. 
This kind of problem of recovering the additive structure of function 
fields has been treated in [Uchida] and [Tamagawa1], using certain 
auxiliary functions --- the 
minimal elements, i.e., functions with a minimal pole. The arguments loc. cit.
work well in the case of a bijection between multiplicative groups, but fail in our case where 
we only have an embedding $\Cal O_{E_X}^{\times} \hookrightarrow \Cal O_{E_Y}^{\times}$ between 
multiplicative groups. In our case, 
instead of using minimal elements, we can recover the additivity by 
using functions whose divisor has a unique pole. 
In our situation the fact that we can only evaluate functions at all 
but finitely many points of $X$ (or, more precisely, 
all points of $X-E_X$) presents an additional 
difficulty which we overcome, roughly speaking, by passing to an infinite algebraic extension of 
the base field, and using ``infinitely many'' auxiliary functions.

In \S 1, we review some 
(mostly known but partly new) 
results which show that various invariants of the curve 
$X$ (among other things, the set of decomposition groups at closed points 
of $X$) can be recovered group-theoretically, 
starting from $\Pi_X$. We also prove the almost injectivity of the above map 
{}from the set of closed points of $X$ to the set of conjugacy 
classes of decomposition groups. 
In \S 2, we review the main results of Mochizuki's theory of
cuspidalizations of \'etale fundamental groups of proper hyperbolic curves, 
which plays an essential role in this paper. 
In \S 3, we prove our main results, assuming the 
results of \S 4. Finally, 
In \S 4, we investigate the problem of recovering the 
additive structure of functions mentioned above. Using the above 
``one pole argument'', we prove 
Proposition 4.4, which is used in \S 3.

\definition{Acknowledgment} 
The authors would like to thank Shinichi Mochizuki very much 
for various useful discussions on his theory of cuspidalizations. 
The first author holds an EPSRC advanced research fellowship
GR/R75861/02, and would like very much to thank EPSRC for its support. 
This work was done during a visit of his to the Research Institute 
for Mathematical Sciences (RIMS) in Kyoto. 
He would like very much to thank the staff members of RIMS 
for their hospitality.
\enddefinition

\subhead
\S 1. Characterization of Decomposition Groups
\endsubhead

Let $X$ be a proper, smooth, geometrically connected 
curve over a finite field $k=k_X$ of characteristic 
$p=p_X>0$. Write $K=K_X$ for the function field 
of $X$. 

Let $S$ be a (possibly empty) finite set of closed 
points of $X$, and set $U=U_S\defeq X-S$. We assume that 
$U$ is hyperbolic. 

Fix a separable closure $K^{\sep}=K^{\sep}_X$ of $K$, 
and write $\overline{k}=\overline{k_X}$ 
for the algebraic closure of $k$ in $K^{\sep}$. 
Write 
$$G_{K}\defeq\Gal(K^{\sep}/K),$$
$$G_{k}\defeq\Gal(\overline{k}/k)$$ 
for the absolute Galois groups of $K$ and $k$, respectively. 

The tame fundamental group $\pi_1^t(U)$ 
with respect to the base point defined by $K^{\sep}$ 
(where ``tame'' is with respect to the complement of $U$ in $X$)
can be naturally identified with a quotient of 
$G_{K}$. 
Write $\Gal(K^t_U/K)$ for 
this quotient. 
(In case $S=\emptyset$, we also write $K^{\ur}_U$ for $K^t_U$.) 
It is easy to see that 
$K^t_U$ contains 
$K\overline{k}$. 

Let $\Sigma=\Sigma_X$ be a set of prime numbers that contains at least
one prime number different from $p$. Write 
$$\Sigma ^{\dag}\defeq \Sigma - \{p\}.$$ 
Thus, $\Sigma ^{\dag}\neq \emptyset$ by our assumption.
Denote by $\hat \Bbb Z^{\Sigma^\dag}$ the maximal pro-$\Sigma ^{\dag}$ 
quotient of $\hat \Bbb Z$. 
Set $\Sigma'=\Sigma'_X=\Primes-\Sigma_X$. We say that 
$\Sigma$ is cofinite if $\sharp(\Sigma')<\infty$. 

We define ${\tilde K}_U$ to be the 
maximal pro-$\Sigma$ 
subextension of $K\overline{k}$ in $K^{t}_U$. 
Now, set 
$$\Pi_U=\Gal({\tilde K}_U/K),$$ 
which is a quotient of 
$\pi_1^t(U)=\Gal(K^{t}_U/K)$. 
This fits into the exact sequence 
$$1\to \Delta_U \to \Pi_U \overset{\pr_U}\to{\to} G_{k}\to 1.$$ 
Here, $\Delta_U$ is the maximal 
pro-$\Sigma$ quotient of 
$\pi_1^t({\overline U})$, where, 
for a $k$-scheme $Z$, we set $\overline Z \defeq Z \times_k \overline k$. 

Define ${\tilde X}_U$ to be the integral closure 
of $X$ in ${\tilde K}_U$. 
Define 
${\tilde U}$ to be the integral closure 
of $U$ in ${\tilde K}_U$, 
which can be naturally identified with the inverse image 
(as an open subscheme) of $U$ in ${\tilde X}_U$. Define 
${\tilde S}_U$ to be the inverse image (as a set) 
of $S$ in ${\tilde X}_U$. 

For a scheme $Z$, write $Z^{\cl}$ for the set of closed points of $Z$. 
Then we have 
$$X^{\cl}=U^{\cl}\coprod S,$$ 
$$({\tilde X}_U)^{\cl}={\tilde U}^{\cl}\coprod {\tilde S}_U.$$
Moreover, 
$({\tilde X}_U)^{\cl}$ admits a 
natural action of $\Pi_U$, and 
the corresponding quotient can be naturally identified with 
$X^{\cl}$. 

For each $\tilde x \in ({\tilde X}_U)^{\cl}$, we define the 
decomposition group $D_{\tilde x}\subset \Pi_U$ 
(respectively, the inertia group $I_{\tilde x}\subset D_{\tilde x}$) 
to be the stabilizer at $\tilde x$ of the 
natural action of $\Pi_U$ on $({\tilde X}_U)^{\cl}$ 
(respectively, the kernel of the natural action of $D_{\tilde x}$ 
on $k({\tilde x})=\overline{k(x)}=\overline{k}$). These 
groups fit into the following commutative diagram in which 
both rows are exact:
$$
\matrix
1 &\to& I_{\tilde x} &\to& D_{\tilde x} &\to& G_{k(x)} &\to& 1\\
&&&&&&&&\\
&& \cap && \cap && \cap && \\
&&&&&&&&\\
1 &\to& \Delta_U &\to& \Pi_U &\to& G_{k} &\to& 1
\endmatrix
$$
Moreover, 
$I_{\tilde x}=\{1\}$ (respectively,  
$I_{\tilde x}$ is (non-canonically) isomorphic to 
${\hat {\Bbb Z}}^{\Sigma^\dag}$), 
if ${\tilde x}\in {\tilde U}^{\cl}$ (respectively, 
${\tilde x}\in {\tilde S}_U$). 
Since $I_{\tilde x}$ is normal in $D_{\tilde x}$, 
$D_{\tilde x}$ acts on $I_{\tilde x}$ by conjugation. 
Since $I_{\tilde x}$ is abelian, this action factors through 
$D_{\tilde x}\to G_{k(x)}$ and induces a natural 
action of $G_{k(x)}$ on $I_{\tilde x}$. 

\proclaim{Lemma 1.1} Assume 
${\tilde x}\in {\tilde S}_U$. Then: 

\noindent
{\rm (i)} The subgroup $I_{\tilde x}^{G_{k(x)}}$ of $I_{\tilde x}$ 
that consists of elements fixed by the $G_{k(x)}$-action is trivial. 

\noindent
{\rm (ii)} Assume moreover that $\Sigma$ is cofinite. Then 
the kernel of the action $G_{k(x)}\to \Aut(I_{\tilde x})$ 
is trivial. 
\endproclaim

\demo{Proof} 
By assumption, $I_{\tilde x}\simeq {\hat {\Bbb Z}}^{\Sigma^{\dag}}$, and it is 
well-known that the map $G_{k(x)}\to \Aut(I_{\tilde x})
=({\hat {\Bbb Z}}^{\Sigma^{\dag}})^{\times}$ coincides with the cyclotomic 
character and sends the $\sharp(k(x))$-th power 
Frobenius element $\varphi_{k(x)}\in G_{k(x)}$, which is a (topological) generator 
of $G_{k(x)}$, to $\sharp(k(x))\in ({\hat {\Bbb Z}}^{\Sigma^{\dag}})^{\times}$. 
The assertion of (i) follows from this, since $\sharp(k(x))-1$ 
is not a zero divisor of the ring ${\hat {\Bbb Z}}^{\Sigma^{\dag}}$. The assertion 
of (ii) also follows from this, together with a theorem of Chevalley 
([Chevalley], Th\'eor\`eme 1). 
\qed
\enddemo

Let $G$ be a profinite group. Then, define $\Sub(G)$ 
(respectively, $\OSub(G)$) 
to be the set of closed (respectively, open) subgroups of $G$. 

By conjugation, $G$ acts on $\Sub(G)$. More generally, 
let $H$ and $K$ be closed subgroups of $G$ such that 
$K$ normalizes $H$. Then, by conjugation, $K$ acts on 
$\Sub(H)$. We denote by $\Sub(H)_K$ the quotient 
$\Sub(H)/K$ by this action. In particular, 
$\Sub(G)_G$ is the set of conjugacy classes of 
closed subgroups of $G$. 

For any closed subgroups $H,K$ of $G$ with $K\subset H$, 
we have 
a natural inclusion 
$\Sub(K)\subset \Sub(H)$, as well as 
a natural map
$\Sub(H)\to\Sub(K)$, $J\mapsto J\cap K$. By using this latter 
natural map, we define
$$\overline{\Sub}(G)\defeq\underset{H\in\OSub(G)}\to{\varinjlim}
\Sub(H).$$
Observe that $\overline{\Sub}(G)$ can be identified with 
the set of commensurate classes of closed subgroups of $G$. 
(Closed subgroups $J_1$ and $J_2$ of $G$ are called commensurate 
(to each other), if $J_1\cap J_1$ is open both in $J_1$ and in $J_2$.) 

With these notations, we obtain natural maps
$$D=D[U]:({\tilde X}_U)^{\cl}\to\Sub(\Pi_U), {\tilde x}\mapsto D_{\tilde x},$$
$$I=I[U]:({\tilde X}_U)^{\cl}\to\Sub(\Delta_U)\subset
\Sub(\Pi_U), {\tilde x}\mapsto I_{\tilde x},$$
which fit into the commutative diagram
$$\CD
({\tilde X}_U)^{\cl} @>{D}>> \Sub(\Pi_U) \\
\| @. @VVV \\
({\tilde X}_U)^{\cl} @>{I}>> \Sub(\Delta_U), 
\endCD
$$
where the vertical arrow stands for the natural map 
$\Sub(\Pi_U)\to \Sub(\Delta_U)$, 
$J\mapsto J\cap \Delta_U$. 
By composition with the natural map 
$\Sub(\Pi_U)\to\overline{\Sub}(\Pi_U)$, $D,I$ yield 
$$\overline{D}=\overline{D}[U]:({\tilde X}_U)^{\cl}\to\overline{\Sub}(\Pi_U),$$
$$\overline{I}=\overline{I}[U]:({\tilde X}_U)^{\cl}\to\overline{\Sub}
(\Delta_U)\subset \overline{\Sub}(\Pi_U).$$

\definition{Remark 1.2} Unlike the case of $D,I$, 
the maps $\overline{D},\overline{I}$ are essentially unchanged 
if we replace $U$ by any covering corresponding to an open 
subgroup of $\Pi_U$. 
\enddefinition

Since 
the maps $D,I$ 
are $\Pi_U$-equivariant, they 
induce natural maps 
$$D_{\Pi_U}=D[U]_{\Pi_U}:X^{\cl}\to\Sub(\Pi_U)_{\Pi_U},$$
$$I_{\Pi_U}=I[U]_{\Pi_U}:X^{\cl}\to\Sub(\Delta_U)_{\Pi_U}
\subset\Sub(\Pi_U)_{\Pi_U},$$
respectively. 

\definition{Definition 1.3} 
Let $f:A\to B$ be a map of sets. 

\noindent
(i) We define $\mu_f: B\to {\Bbb Z}\cup\{\infty\}$ by 
$\mu_f(b)=\sharp(f^{-1}(b))$. (Thus, 
$f$ is injective (respectively, surjective) if 
$\mu_f(b)\leq 1$ (respectively, $\mu_f(b)\geq 1$) for 
any $b\in B$. We also have  
$f(A)=\{b\in B\mid \mu_f(b)\geq 1\}$.)

\noindent
(ii) We say that $f$ is quasi-finite, 
if $\mu_f(b)<\infty$ for any $b\in B$. 

\noindent
(iii) We say that an element $a$ of $A$ is an exceptional element 
of $f$ (in $A$), if $\mu_f(f(a))>1$. We refer to the set of 
exceptional elements as the exceptional set. 

\noindent
(iv) We say that a pair $(a_1,a_2)$ of elements of $A$ is an exceptional 
pair of $f$ (in $A$), if $a_1\neq a_2$ and $f(a_1)=f(a_2)$ hold. 

\noindent
(v) We say that $f$ is almost injective (in the strong sense), 
if the exceptional set of $f$ is finite. 
(Observe that almost injectivity implies quasi-finiteness.) 
\enddefinition

\proclaim{Lemma 1.4}
Let $f:A\to B$ and $g:B\to C$ be maps of sets. Then we have: 
$$
\matrix
\text{Both $f$ and $g$ are quasi-finite (respectively, almost injective).}\\
\Downarrow\\
\text{$g\circ f$ is quasi-finite (respectively, almost injective).}\\
\Downarrow\\
\text{$f$ is quasi-finite (respectively, almost injective).}
\endmatrix
$$
\endproclaim

\demo{Proof} 
Easy. 
\qed
\enddemo

\definition{Definition 1.5}
Denote by $E_{\tilde U}$ the exceptional set 
of $\overline{D}$ in $({\tilde X}_U)^{\cl}$ 
\enddefinition

\definition{Definition 1.6}
Let $G$ be a profinite group and $H$ a closed subgroup. Then 
we denote by $Z_G(H)$, $N_G(H)$ and $C_G(H)$ the 
centralizer, the normalizer and the commensurator, respectively, 
of $H$ in $G$. Namely, 
$$
\matrix
\format\c&\l\\
Z_G(H)&=\{g\in G \mid \text{$ghg^{-1} = h$ for any $h\in H$}\},\\
\cap\\
N_G(H)&=\{g\in G \mid gHg^{-1} = H\}\supset H,\\
\cap \\
C_G(H)&=\{g \in G \mid \text{$gHg^{-1}$ and $H$ are commensurate}\}.
\endmatrix$$
\enddefinition

\proclaim{Lemma 1.7} 
Let $Z$ be a closed subgroup of $\Pi_U$ such that $\pr_U(Z)$ is open 
in $G_k$ and that $\pr_U|_Z$ is injective. Then 
$\pr_U$ induces an injection $C_{\Pi_U}(Z) \hookrightarrow G_{k}$,  
and we have $C_{\Pi_U}(Z)=N_{\Pi_U}(Z)=Z_{\Pi_U}(Z)
\supset Z$
and 
$(C_{\Pi_U}(Z): Z)<\infty$. 
\endproclaim

\demo{Proof} 
Take any $\sigma\in C_{\Pi_U}(Z)\cap \Delta_U$. 
Thus, $Z_0\defeq Z\cap \sigma Z \sigma^{-1}$ 
is open both in $Z$ and in $\sigma Z\sigma^{-1}$. 
We claim that $\sigma$ commutes with any element $\tau$ of $Z_0$. 
Indeed, first, observe that 
$\tau\in Z_0\subset \sigma Z\sigma^{-1}$
and 
$\sigma\tau\sigma^{-1}\in 
\sigma Z_0\sigma^{-1}\subset \sigma Z\sigma^{-1}$ 
hold. Or, equivalently, 
$\sigma^{-1}\tau\sigma, \tau \in Z$. 
Second, observe that 
$\pr_U(\sigma^{-1}\tau\sigma)=\pr_U(\tau)$
holds, since $\pr_U(\sigma)=1$. Now
since $\pr_U|_Z$ is injective, the equality 
$\pr_U(\sigma^{-1}\tau\sigma)=\pr_U(\tau)$
implies $\sigma^{-1}\tau\sigma=\tau$, as desired. 

Next, we prove $\sigma=1$. To see this, suppose 
$\sigma\neq 1$ and take any 
sufficiently small, open characteristic subgroup 
${\overline N}$ of $\Delta_U$ such that $\sigma\not\in {\overline N}$. 
Set ${\overline H}\defeq \langle {\overline N}, \sigma\rangle 
\subset \Delta_U$. Then the image of $\sigma$ in 
${\overline H}^{\ab}$ is nontrivial. (Indeed, the image of 
$\sigma$ in ${\overline H}/{\overline N}$ is nontrivial by definition. 
Since ${\overline H}/{\overline N}$ is cyclic, the 
surjection ${\overline H}\to {\overline H}/{\overline N}$ 
factors through the surjection ${\overline H}\to{\overline H}^{\ab}$.) 
Observe that $Z_0$ normalizes ${\overline H}$, since 
${\overline N}$ is characteristic in $\Delta_U$ and 
$\sigma$ commutes with $Z_0$. So, the open subgroup 
$H\defeq\langle {\overline H}, Z_0\rangle$ of $\Pi_U$ 
can be regarded as a semidirect-product extension 
of $Z_0$ by ${\overline H}$ and satisfies 
$H\cap\Delta_U={\overline H}$. Now, the image 
of $\sigma$ in ${\overline H}^{\ab}$ is nontrivial and 
fixed by the action of $Z_0$. This is impossible, 
as can be easily seen by observing the 
Frobenius weights in the action of $Z_0$, 
or of $\pr_U(Z_0)$, which is an open subgroup 
of $G_{k}$. 

Thus, we have proved $\sigma=1$, and the first assertion follows from this. 
In particular, $C_{\Pi_U}(Z)$ ($\hookrightarrow G_k$) is abelian, hence 
the second assertion follows. Finally, since $\pr_U$ induces an isomorphism 
$C_{\Pi_U}(Z)\overset{\sim}\to{\to}\pr_U(C_{\Pi_U}(Z))$ and 
$\pr_U(Z)$ is open in $G_k$, the third assertion holds. 
\qed\enddemo

The first main result in this \S\  is: 

\proclaim{Proposition 1.8} {\rm (i)} 
$\overline{I}|_{{\tilde S}_U}:{\tilde S}_U\to\overline{\Sub}(\Pi_U)$ is injective. 

\noindent
{\rm (ii)} 
$E_{\tilde U}$ is disjoint from ${\tilde S}_U$. 
(Or, equivalently, $E_{\tilde U}\subset {\tilde U}^{\cl}$.)

\noindent
{\rm (iii)} 
Let $\overline \rho$ denote the natural morphism 
${\tilde X}_U\to {\overline X}$. Then, 
for each ${\overline x}\in {\overline X}^{\cl}$, 
$\overline{D}|_{{\overline \rho}^{-1}({\overline x})}$ is injective. 

\noindent
{\rm (iv)} Let $\rho$ denote the natural morphism 
${\tilde X}_U\to X$. Then, 
for each $x\in X^{\cl}$, 
$\overline{D}|_{\rho^{-1}(x)}$ is quasi-finite. 
If, moreover, $k(x)=k$ holds (i.e., $x$ is a $k$-rational 
point of $X$), $\overline{D}|_{\rho^{-1}(x)}$ is injective. 

\noindent
{\rm (v)} $E_{\tilde U}$ is $\Pi_U$-stable. 

Assume, moreover, that $\Sigma$ is cofinite. Then: 

\noindent
{\rm (vi)} The quotient $E_{\tilde U}/\Pi_U$ is finite. 

\noindent
{\rm (vii)} $\overline{D}:({\tilde X}_U)^{\cl}\to\overline{\Sub}(\Pi_U)$ 
is quasi-finite. 
\endproclaim

\demo{Proof} 
\noindent
(i) Take any ${\tilde x}, {\tilde x}' \in {\tilde S}_U$, and assume 
${\tilde x}\neq {\tilde x}'$. Then there exists an open subgroup 
$H_0$ of $\Pi_U$, such that the following holds: Let $U_0$ denote the 
covering of $U$ corresponding to $H_0\subset\Pi_U$ 
and $X_0$ the integral closure of $X$ in $U_0$ (i.e., $X_0$ is the smooth 
compactification of $U_0$), then the images $x_0,x'_0$ of 
${\tilde x}, {\tilde x}'$ in $X_0$ are distinct from each other. Moreover, 
by replacing $H_0$ by a smaller open subgroup if necessary, 
we may assume that 
the cardinality of the point set ${\overline X}_0-{\overline U}_0$ is $\geq 3$ 
(see, e.g., [Tamagawa1], Lemma (1.10)). 

Now, to show the desired injectivity, 
it suffices to prove that
$I_{\tilde x}\cap H_1 \neq I_{{\tilde x}'}\cap H_1$ holds 
for any open subgroup $H_1$ of $H_0$. 
Let $U_1$ denote the 
covering of $U$ corresponding to $H_1\subset\Pi_U$ 
and $X_1$ the integral closure of $X$ in $U_1$. Then, by the 
choice of $H_1$, we see that the images 
of ${\tilde x}, {\tilde x}'$ in 
${\overline S}_1\defeq{\overline X}_1-{\overline U}_1$ are distinct from each other 
and that the cardinality of ${\overline S}_1$ is $\geq 3$. Then 
it is easy to see that the images of 
$I_{\tilde x}\cap H_1,I_{{\tilde x}'}\cap H_1$ 
in ${\overline H}_1^{\ab}$ are isomorphic to ${\hat {\Bbb Z}}^{\Sigma^{\dag}}$ 
and that the intersection of these images is $\{0\}$. 
(Observe (the pro-$\Sigma^{\dag}$ part of) 
exact sequence (1-5) in [Tamagawa1].) Thus, 
a fortiori, $I_{\tilde x}\cap H_1 \neq I_{{\tilde x}'}\cap H_1$ holds, 
as desired. 

\noindent
(ii) Take any ${\tilde x}\in {\tilde S}_U$ 
and ${\tilde x}'\in({\tilde X}_U)^{\cl}$, 
such that ${\tilde x}\neq {\tilde x}'$ holds, 
then we shall prove that the images of
${\tilde x},{\tilde x}'$ by $\overline{D}$ 
are distinct from each other. To see this, 
it suffices, by definition, to prove that, for any open subgroup 
$H$ of $\Pi_U$, the images 
$D_{\tilde x}\cap H,D_{{\tilde x}'}\cap H$
of ${\tilde x},{\tilde x}'$ 
in $\Sub(H)$ are distinct from each other. Now, replacing 
$U$ by the covering of $U$ corresponding to $H\subset \Pi_U$, 
it suffices to prove that 
$D_{\tilde x},D_{{\tilde x}'}$ are distinct from each other. 
Now, recall that
$D_{\tilde x}\cap \Delta_U=I_{\tilde x},\  
D_{{\tilde x}'}\cap \Delta_U=I_{{\tilde x}'}.$
Thus, if ${\tilde x}'\in {\tilde S}_U$, 
the last assertion follows from (i). On the other hand, 
if ${\tilde x}'\in {\tilde U}^{\cl}$, the last assertion 
follows from the fact 
$I_{\tilde x}\simeq 
{\hat {\Bbb Z}}^{\Sigma^{\dag}}$, $I_{{\tilde x}'}=\{1\}$. 

\noindent
(iii) If ${\overline x}\in {\overline S}\defeq{\overline X}-{\overline U}$, 
the assertion 
follows from (ii). So, we may and shall 
assume ${\overline x}\in {\overline U}^{\cl}$. 
Take any ${\tilde x}, {\tilde x}' \in {\overline \rho}^{-1}({\overline x})$. 
Then there exists $\sigma\in\Delta_U$ such that 
${\tilde x}'=\sigma{\tilde x}$ holds. 
(Such $\sigma$ is unique by the assumption 
${\overline x}\in {\overline U}^{\cl}$, though we do not use this fact 
in the proof.) 
Now, suppose that the images of ${\tilde x}, {\tilde x}'$ 
by $\overline{D}$ coincide with each other. Namely, 
$D_{\tilde x}$ and 
$D_{{\tilde x}'}=D_{\sigma{\tilde x}}=\sigma D_{\tilde x}\sigma^{-1}$ 
are commensurate to each other. Thus, $\sigma\in C_{\Pi_U}(D_{\tilde x}) 
\cap \Delta_U$, and it follows from Lemma 1.7 that $\sigma=1$ holds, hence 
${\tilde x}'=\sigma{\tilde x}={\tilde x}$. 
Namely, $\overline{D}|_{{\overline \rho}^{-1}({\overline x})}$ is injective, 
as desired. 

\noindent
(iv) Let $\pi$ denote the natural morphism ${\overline X}\to X$, 
so that $\rho=\pi\circ {\overline \rho}$ holds. Since 
$\sharp(\pi^{-1}(x))=[k(x):k]<\infty$, the assertions follow 
directly from (iii). 

\noindent
(v) This follows from the fact that $\overline{D}$ is 
$\Pi_U$-equivariant. 

\noindent
(vi) To prove this (assuming that $\Sigma$ is cofinite),  we may replace 
$U$ by any covering corresponding to an open subgroup of $\Pi_U$. 
Thus, we may assume that the genus of $X$ is $>1$ and that 
$X$ is non-hyperelliptic. (See, e.g., [Tamagawa1], Lemma (1.10) 
for the former, 
and either [Tamagawa3], \S 2 or the proof (in characteristic zero) 
of [Mochizuki], Lemma 10.4(4) for the latter.) We shall prove 
that $\rho(E_{\tilde U})$, which can be identified with 
$E_{\tilde U}/\Pi_U$, is finite, or, more strongly, that 
${\overline \rho}(E_{\tilde U})$, which can be identified with 
$E_{\tilde U}/\Delta_U$, is finite. 

Take any pair of elements 
${\tilde x}, {\tilde x}'\in {\tilde U}^{\cl}$, 
and denote by ${\overline x},{\overline x}'$ the images of 
${\tilde x}, {\tilde x}'$ in ${\overline U}^{\cl}$, respectively. 
The images $\pr_U(D_{\tilde x})$ and $\pr_U(D_{{\tilde x}'})$ 
are open in $G_{k}$, hence so is the 
intersection $G_0\defeq
\pr_U(D_{\tilde x})\cap\pr_U(D_{{\tilde x}'})$. 
Let $s,s'$ be the inverse maps of the isomorphisms 
$\pr_U|_{D_{\tilde x}}: D_{\tilde x}\to \pr_U(D_{\tilde x})$, 
$\pr_U|_{D_{{\tilde x}'}}: D_{{\tilde x}'}\to \pr_U(D_{{\tilde x}'})$, 
respectively. Then, it is well-known and easy to see that 
the map $\phi: G_0\to\Delta_U$, 
$\gamma\mapsto s(\gamma)s'(\gamma)^{-1}$ 
is a continuous $1$-cocycle (with respect to 
the left, conjugacy action of $G_0$ on $\Delta_U$ via 
the section $s'$). Thus, $\phi$ defines a cohomology class 
in $H^1(G_0, \Delta_U)$. We denote by 
$\phi_{0,X}^{\ab}=\phi_{0,X}^{\ab}(\tilde x,\tilde x')$ 
the image of 
this class in $H^1(G_0, \Delta_{X}^{\ab})$. 
(Note that the $G_0$-action on 
$\Delta_{X}^{\ab}$ induced by that on $\Delta_U$ 
extends to a canonical $G_{k}$-action, 
hence, in particular, is independent of the choice 
of ${\tilde x}'$.) Moreover, we set 
$${\Cal H}_X\defeq 
\varinjlim_{G\in\OSub(G_k)}
H^1(G, \Delta_{X}^{\ab})$$
(where the transition maps are restriction maps)
and denote by $\phi_X^{\ab}=\phi_X^{\ab}({\tilde x},{\tilde x}')$ 
the image of $\phi_{0,X}^{\ab}$ in ${\Cal H}_X$. 

On the other hand, it is well-known that $\Delta_{X}^{\ab}$ 
is canonically isomorphic as a $G_{k}$-module to 
the pro-$\Sigma$ part $T(J)^{\Sigma}$ 
of the full Tate module $T(J)$ 
of the Jacobian variety $J$ (tensored with $\overline{k}$) 
of $X$. Thus, by Kummer theory (for the abelian variety $J$), 
we obtain an injective map 
$J(k_G)/(J(k_G)\{\Sigma'\})\to H^1(G, \Delta_{X}^{\ab})$, 
where $G$ is an open subgroup of $G_{k}$, $k_G$ is 
the finite extension of $k$ corresponding to $G$, 
and, for an abelian group $M$, $M\{\Sigma'\}$ stands 
for the subgroup of torsion elements $a$ of $M$ such 
that every prime divisor of the order of $a$ 
belongs to $\Sigma'$. 
(In fact, the above injective map is bijective by Lang's theorem, 
though we do not use this fact in the proof.) By taking 
the inductive limit, we obtain an injective map 
$J(\overline{k})/ (J(\overline{k})\{\Sigma'\})\to 
{\Cal H}_X$. 
Now, it is widely known that the image in ${\Cal H}_X$ 
of the class of ${\overline x}-{\overline x}'$ 
in ${\overline J}^{\cl}=J(\overline{k})$ 
coincides with $\phi_X^{\ab}$. For this, see [NT], Lemma (4.14). 
(See also [Nakamura2], 2.2 and [Tamagawa1], Lemma (2.6).) 

Suppose moreover that $(\tilde x,{\tilde x}')$ is an exceptional 
pair of $\overline{D}$. Then it follows from the various definitions 
that $\phi_X^{\ab}\in{\Cal H}_X$ is trivial. Therefore the class 
of ${\overline x}-{\overline x}'$ in 
$J(\overline{k})/ (J(\overline{k})\{\Sigma'\})$ 
is trivial, or, equivalently, the class 
$\cl({\overline x}-{\overline x}')$ in $J(\overline{k})$ falls 
in $J(\overline{k})\{\Sigma'\}$. On the other hand, 
by (iii), ${\overline x}\neq{\overline x}'$ hold, or, equivalently 
(by the assumption that the genus of $X$ is $>1$), 
$\cl({\overline x}-{\overline x}')\neq 0$. 

Consider the morphism $\delta: X\times X\to J$, 
$(P,Q)\mapsto\cl(P-Q)$. We claim: 

\definition{Claim 1.9} 
(i) $\delta|_{X\times X-\iota(X)}$ is injective 
(on $\overline{k}$-valued points), where 
$\iota :X\to X\times X$ is the diagonal morphism. 

\noindent
(ii) The image $W$ of $\delta$ does not contain any translate 
of positive-dimensional abelian subvariety of $J$. 
\enddefinition

Indeed, for (i), suppose that $(P,Q),(P',Q')\in 
(X\times X-\iota(X))(\overline{k})$ have the same image 
under $\delta$. Namely, the divisors $P-Q$ and $P'-Q'$ 
are linearly equivalent: $P-Q\sim P'-Q'$, or, equivalently, 
$P+Q'\sim P'+Q$. Since we have assumed that $X$ is of genus 
$>1$ and non-hyperelliptic, this implies that 
$P+Q'$ and $P'+Q$ coincide with each other as divisors. 
This implies that either $P=P',Q=Q'$ or $P=Q,P'=Q'$ holds. 
The former implies that $(P,Q)=(P',Q')$, as desired, and 
the latter implies that $(P,Q),(P',Q')\in\iota(X)$, which 
contradicts the assumption. For (ii), suppose that 
$W$ contains a translate $B'$ of some positive-dimensional 
abelian subvariety $B$ of $J$. As $\dim(W)\leq\dim(X\times X)=2$, 
we have $\dim(B')\leq 2$, i.e., either $\dim(B')=2$ or 
$\dim(B')=1$. The former implies that $B'=W$, since $W$ is 
irreducible of dimension $\leq 2$ as the image of $X\times X$. 
Since $0\in W=B'$, we conclude $W=B'=B$. Now, since $J$ is generated 
by $W$, we must have $J=B$. This implies that the genus of $X$ 
(i.e., $\dim(J)$) is $2$, which implies that $X$ is hyperelliptic. 
This contradicts the assumption. So, suppose $\dim(B')=1$. 
By (i), we see that $\delta$ induces a bijective morphism 
$X\times X-\Delta(X)\to W-\{0\}$. From this, we deduce that 
there exists a finite radicial covering $B''$ of $B'$ that 
admits a non-constant morphism to $X\times X$. 
In particular, considering a suitable one of two projections, 
we see that $B''$ admits a non-constant morphism to $X$. 
This is absurd, since the genus of $B''$ (respectively, $X$) is 
$1$ (respectively, $>1$). This completes the proof of Claim 1.9. 

By Claim 1.9(ii) and [Boxall] (which is the most nontrivial 
ingredient of the proof of Proposition 1.8), we see that 
$W(\overline{k})\cap (J(\overline{k})\{\Sigma'\})$ is finite. 
Now, by Claim 1.9(i), we conclude that there exists a finite 
subset $\Cal S$ of $(X\times X)(\overline{k})$ that contains 
any pair $({\overline x},{\overline x}')$ as above. This implies the 
desired assertion that 
${\overline \rho}(E_{\tilde U})$ is a finite set. 

\noindent
(vii) Note that $\rho(E_{\tilde U})$ 
can be identified with $E_{\tilde U}/\Pi_U$ by 
(v). Thus, the assertion of (vii) directly follows from 
(vi) and the first part of (iv). 
\qed
\enddemo

\definition{Definition 1.10}
We define $E_U$ to be the image of $E_{\tilde U}$ in $X^{\cl}$. 
(This can be identified with 
$E_{\tilde U}/\Pi_U$. Thus, if $\Sigma$ is cofinite, then it 
is finite by Proposition 1.8(vi).) 
\enddefinition

\proclaim{Corollary 1.11} 
{\rm (i)} $D_{\Pi_U}|_{X^{\cl}-E_U}:X^{\cl}-E_U\to\Sub(\Pi_U)_{\Pi_U}$ 
is injective. 

\noindent
{\rm (ii)} $E_U$ is disjoint from $S$. 
(Or, equivalently, $E_U\subset U^{\cl}$.) 

Assume, moreover, that $\Sigma$ is cofinite. Then: 

\noindent
{\rm (iii)} $D_{\Pi_U}:X^{\cl}\to\Sub(\Pi_U)_{\Pi_U}$ is almost injective. 
\endproclaim

\demo{Proof} 
(i) As $D|_{({\tilde X}_U)^{\cl}-E_{\tilde U}}: 
({\tilde X}_U)^{\cl}-E_{\tilde U}\to \Sub(\Pi_U)$ is 
injective by definition and $\Pi_U$-equivariant, 
its quotient by $\Pi_U$, which is naturally identified with 
$D_{\Pi_U}|_{X^{\cl}-E_U}:X^{\cl}-E_U\to\Sub(\Pi_U)_{\Pi_U}$, 
is also injective. This completes the proof. 

\noindent
(ii) This follows from Proposition 1.8(ii). 

\noindent
(iii) This follows from (i) and the fact that $E_U$ is finite 
(Proposition 1.8(vi)). 
\qed
\enddemo

\proclaim{Corollary 1.12} 
{\rm (i)} For each $\tilde x\in {\tilde U}^{\cl}$, 
$\pr_U$ induces an injection 
$$C_{\Pi_U}(D_{\tilde x}) \hookrightarrow G_{k},$$ 
and we have 
$$C_{\Pi_U}(D_{\tilde x})=N_{\Pi_U}(D_{\tilde x})=Z_{\Pi_U}(D_{\tilde x})
\supset D_{\tilde x}$$
and 
$$(C_{\Pi_U}(D_{\tilde x}): D_{\tilde x})<\infty.$$

If, moreover, ${\tilde x}\in {\tilde U}^{\cl}-E_{\tilde U}$, 
we have 
$$C_{\Pi_U}(D_{\tilde x})=N_{\Pi_U}(D_{\tilde x})=Z_{\Pi_U}(D_{\tilde x})
= D_{\tilde x}.$$

\noindent
{\rm (ii)} For each $\tilde x\in {\tilde S}_U$, we have 
$$C_{\Pi_U}(D_{\tilde x})=N_{\Pi_U}(D_{\tilde x})=D_{\tilde x},\  
Z_{\Pi_U}(D_{\tilde x})=Z_{D_{\tilde x}}(D_{\tilde x})$$
and 
$$C_{\Pi_U}(I_{\tilde x})=N_{\Pi_U}(I_{\tilde x})=D_{\tilde x},\  
Z_{\Pi_U}(I_{\tilde x})=Z_{D_{\tilde x}}(I_{\tilde x}).$$
(Note that, if $\Sigma$ is cofinite, then 
$Z_{D_{\tilde x}}(D_{\tilde x})=1$ and 
$Z_{D_{\tilde x}}(I_{\tilde x})=I_{\tilde x}$ 
by Lemma 1.1.) 

\noindent
{\rm (iii)} Assume, moreover, that $\Sigma$ is cofinite. Then 
there exists an open subgroup $G_0$ of $G_{k}$, such that,  
for 
any open subgroup $H$ of $\pr_U^{-1}(G_0)$ and 
any element 
$\tilde x$ of $({\tilde X}_U)^{\cl}={\tilde U}^{\cl}\coprod {\tilde S}_U$, 
we have 
$$C_{H}(D_{\tilde x}\cap H)=N_{H}(D_{\tilde x}\cap H)=
D_{\tilde x}\cap H,$$
$$Z_{H}(D_{\tilde x}\cap H)=
\cases
D_{\tilde x}\cap H, &\text{for ${\tilde x}\in {\tilde U}^{\cl}$},\\
\{1\}, &\text{for ${\tilde x}\in {\tilde S}_U$}.
\endcases
$$
In other words, if we replace $U$ by a covering corresponding to 
such $H$, we have, for any ${\tilde x}\in ({\tilde X}_U)^{\cl}$, 
$$C_{\Pi_U}(D_{\tilde x})=N_{\Pi_U}(D_{\tilde x})=
D_{\tilde x},$$
$$Z_{\Pi_U}(D_{\tilde x})=
\cases
D_{\tilde x}, &\text{for ${\tilde x}\in {\tilde U}^{\cl}$},\\
\{1\}, &\text{for ${\tilde x}\in {\tilde S}_U$}.
\endcases
$$
\endproclaim

\demo{Proof} 
First, 
since 
$\overline{D}|_{({\tilde X}_U)^{\cl}-E_{\tilde U}}: 
({\tilde X}_U)^{\cl}-E_{\tilde U}\to \overline{\Sub}(\Pi_U)$ 
is injective and $\Pi_U$-equivariant, 
we see that $C_{\Pi_U}(D_{\tilde x})=D_{\tilde x}$ holds 
for any ${\tilde x}\in({\tilde X}_U)^{\cl}-E_{\tilde U}$. 

\noindent 
(i) 
The first assertion follows from Lemma 1.7. The second assertion 
follows from the first assertion and 
the fact shown at the beginning of the proof. 

\noindent
(ii) Let ${\tilde x}\in {\tilde S}_U$. Then, ${\tilde x}\not\in E_{\tilde U}$ 
by Proposition 1.8(ii). Thus, we have 
$C_{\Pi_U}(D_{\tilde x})=N_{\Pi_U}(D_{\tilde x})=D_{\tilde x}$. From this, 
we also have $Z_{\Pi_U}(D_{\tilde x})= Z_{D_{\tilde x}}(D_{\tilde x})$. 

Next, by Proposition 1.8(i), the map 
$\overline{I}|_{{\tilde S}_U}: {\tilde S}_U\to \overline{\Sub}(\Pi_U)$ 
is injective. Since this map is also $\Pi_U$-equivariant, 
we see that $C_{\Pi_U}(I_{\tilde x})=D_{\tilde x}$. As 
$C_{\Pi_U}(I_{\tilde x})\supset N_{\Pi_U}(I_{\tilde x})
\supset D_{\tilde x}$, we have 
$N_{\Pi_U}(I_{\tilde x})=D_{\tilde x}$. From this, 
we also have 
$Z_{\Pi_U}(I_{\tilde x})= Z_{D_{\tilde x}}(I_{\tilde x})$. But this 
last group coincides with $I_{\tilde x}$ by Lemma 1.1(ii). 

\noindent
(iii) Define $G_0$ to be the intersection (in $G_{k}$) of $G_{k(x)}$ 
for $x\in E_U$. Since $E_U$ is finite by Proposition 1.8(vi), 
$G_0$ is an open subgroup of $G_{k}$. By (i) and (ii), it is easy to see 
that this $G_0$ satisfies the desired properties. 
\qed
\enddemo

%
%
%
%
%

Next, we shall show that various 
invariants and structures of $U$ can be recovered 
group-theoretically (or $\varphi$-group-theoretically) from $\Pi_U$, 
in the following sense. 

\definition{Definition 1.13}
(i) We say that $\Pi=(\Pi, \Delta, \varphi_\Pi)$ is a 
$\varphi$-(profinite )group, 
if $\Pi$ is a profinite group, $\Delta$ is a closed normal subgroup 
of $\Pi$ and $\varphi_\Pi$ is an element of $\Pi/\Delta$. 

\noindent
(ii) An isomorphism from a $\varphi$-group $\Pi=(\Pi, \Delta, \varphi_\Pi)$ 
to another $\varphi$-group $\Pi'=(\Pi', \Delta', \varphi_{\Pi'})$ is 
an isomorphism $\Pi\overset{\sim}\to{\to}\Pi'$ as profinite groups that 
induces $\Delta \overset{\sim}\to{\to} \Delta'$, 
hence also $\Pi/\Delta \overset{\sim}\to{\to} \Pi'/\Delta'$, such that 
the last isomorphism sends $\varphi_\Pi$ to $\varphi_{\Pi'}$. 
\enddefinition

{}From now on, we regard $\Pi_U$ as a $\varphi$-group by 
$\Pi_U=(\Pi_U, \Delta_U,\varphi_k)$, 
where $\varphi_k$ stands for the $\sharp(k)$-th power Frobenius element 
in $G_k=\Pi_U/\Delta_U$. We shall say that an isomorphism 
$\alpha: \Pi_U\isom \Pi_{U'}$ as profinite groups 
is Frobenius-preserving, if $\alpha$ is an isomorphism 
as $\varphi$-groups. 

\definition{Definition 1.14}
(i) Given an invariant $F(U)$ that depends on the isomorphism class 
(as a scheme) of a hyperbolic curve $U$ over a finite field, we say 
that $F(U)$ can be recovered group-theoretically (respectively, 
$\varphi$-group-theoretically) from $\Pi_U$, 
if any isomorphism (respectively, any Frobenius-preserving isomorphism) 
$\Pi_U\isom \Pi_V$ 
implies $F(U)=F(V)$ for two such curves $U,V$. 

\noindent
(ii) Given an additional structure ${\Cal F}(U)$ (e.g., 
a family of subgroups, quotients, elements, etc.) 
on the profinite group $\Pi_U$ that depends 
functorially on a hyperbolic curve $U$ over a finite field 
(in the sense that, for any isomorphism (as schemes) 
between two such curves $U,V$, 
any isomorphism $\Pi_U\isom\Pi_V$ 
induced by this isomorphism 
$U\isom V$ (unique up to composition with inner automorphisms) 
preserves the structures ${\Cal F}(U)$ and ${\Cal F}(V)$), 
we say that ${\Cal F}(U)$ can be recovered group-theoretically 
(respectively, $\varphi$-group-theoretically) from $\Pi_U$, 
if any isomorphism (respectively, any Frobenius-preserving isomorphism) 
$\Pi_U\isom\Pi_V$ 
between two such curves 
$U,V$ preserves the structures ${\Cal F}(U)$ and ${\Cal F}(V)$. 
\enddefinition

\proclaim{Proposition 1.15}
{\rm I.} 
The following invariants and structures can be 
recovered group-theoretically from $\Pi_U$: 

\noindent
{\rm (i)} The subgroup $\Delta_U$ of $\Pi_U$, hence 
the quotient $G_k=\Pi_U/\Delta_U$. 

\noindent
{\rm (ii)} The subsets $\Sigma$ and $\Sigma^\dag$ of $\Primes$. 
\noindent

\smallskip\noindent
{\rm II.} 
The following invariants and structures can be 
recovered $\varphi$-group-theoretically from $\Pi_U$: 

\noindent
{\rm (iii)} The prime number $p$. 

\noindent
{\rm (iv)} The cardinality $\sharp(k)$ 
(or, equivalently, the isomorphism class of the finite field $k$). 

\noindent
{\rm (v)} The genus $g=g_X$ of $X$ and the cardinality $r=r_U\defeq
\sharp({\overline S})$, where ${\overline S}\defeq{\overline X}-{\overline U}$. 

\noindent
{\rm (vi)} The kernel $I_U$ of the natural surjection $\Pi_U\to\Pi_X$ 
(which coincides with the kernel of the natural surjection 
$\Delta_U\to\Delta_{X}$), hence the quotients 
$\Pi_X=\Pi_U/I_U$, $\Delta_{X}=\Delta_U/I_U$. 

\noindent
{\rm (vii)} The cardinalities $\sharp(X(k))$, $\sharp(U(k))$ and 
$\sharp(S(k))$. 

\smallskip\noindent
{\rm III.} Assume, moreover, that $\Sigma$ is cofinite. Then 
the following structure (hence also {\rm (iii)-(vii)} above) can be 
recovered group-theoretically from $\Pi_U$: 

\noindent
{\rm (viii)} The $\sharp(k)$-th power Frobenius element 
$\varphi_{k}\in G_k$. 

\endproclaim

\demo{Proof} 
(i) Similar to [Tamagawa1], Proposition (3.3)(ii). 
(See also [Mochizuki2], Theorem 1.16(ii).) 

\noindent
(ii) Note that $\Delta_U^{\ab}$ is isomorphic to 
$({\hat {\Bbb Z}}^{\Sigma^\dag})^{2g+r+b-1}\times {\Bbb Z}_p^{c}$, 
where $b=b_U$ stands for the second Betti number of $U$, 
i.e., $b=1$ (respectively, $0$) if $r=0$ (respectively, $r>0$), and $c$ stands for 
the $p$-rank of the Jacobian variety of $X$ (respectively, $0$) if 
$p\in\Sigma$ (respectively, $p\not\in\Sigma$). 
(See, e.g., [Tamagawa1], Corollary (1.2).) 
Here, we always have $2g+r+b-1>c\geq 0$. If, moreover, $p\in\Sigma$ 
and if we replace $\Pi_U$ by a suitable open subgroup, then 
we have $c>0$. 
(See, e.g., [Tamagawa1], Lemma (1.9). See also [Tamagawa2], 
Remark (3.11).) From these, it is easy to see that $\Sigma$ and 
$\Sigma^\dag$ can be recovered group-theoretically from $\Pi_U$. 
(See also [Mochizuki2], Theorem 1.16(i).) 

\noindent 
(iii) By conjugation, 
$G_{k}=\Pi_U/\Delta_U$ acts on $(\Delta_U^{\ab})^{\Sigma^\dag}$, 
hence on the 
$\rank_{{\hat {\Bbb Z}}^{\Sigma^\dag}}((\Delta_U^{\ab})^{\Sigma^\dag})$-th 
exterior power $\bigwedge_{{\hat {\Bbb Z}}^{\Sigma^\dag}}^{\max}
(\Delta_U^{\ab})^{\Sigma^\dag}$. 
Thus, we obtain (purely group-theoretically) the character 
$$\rho^{\det}: G_{k}\to 
\Aut(\bigwedge_{{\hat {\Bbb Z}}^{\Sigma^\dag}}^{\max}(\Delta_U^{\ab})^{\Sigma^\dag})
=({\hat {\Bbb Z}}^{\Sigma^\dag})^{\times}.$$ 
As in the proof of [Tamagawa1], Proposition (3.4)(i), we have 
$\rho^{\det}=\epsilon\chi_{\Sigma^\dag}^{g+n+b-1}$, where $\chi_{\Sigma^\dag}$ 
is the pro-$\Sigma^\dag$ cyclotomic character and $\epsilon$ is a certain 
character (depending on $U$) with values in $\{\pm 1\}$. 
Now, $p$ can be characterized by $\rho^{\det}(\varphi_k)\in \pm p^{\Bbb Z}\  
(\subset ({\hat {\Bbb Z}}^{\Sigma^\dag})^{\times})$. 
(See also [Mochizuki2], Remark 1.18.3.) 

\noindent 
(iv) Similar to [Tamagawa1], Proposition (3.4)(iii). 
(See also [Mochizuki2], Remark 1.18.3.) 

\noindent 
(v) Similar to [Tamagawa1], Proposition (3.5). 
(See also [Mochizuki2], Theorem 1.16(i).) 

\noindent 
(vi) Similar to [Tamagawa1], Proposition (3.7). 

\noindent
(vii) Similar to [Tamagawa1], Proposition (3.8). 
More precisely, by the Lefschetz trace formula, 
we have, for any prime $l\in\Sigma^\dag$, 
$$\sharp(X(k))
=\sum_{i=0}^2(-1)^{i}\tr(\varphi_{k}^{-1}\mid H^i_{\et}({\overline X}, {\Bbb Q}_l))
$$
$$
=
\cases
1+\sharp(k), &g=0,\\
\displaystyle{\sum_{i=0}^2}(-1)^{i}\tr(\varphi_{k}^{-1}\mid H^i(\Delta_{X}, {\Bbb Q}_l)), 
&g>0,
\endcases
$$
and 
$$
\align
\sharp(U(k))
&=\sum_{i=0}^2(-1)^{i}\tr(\varphi_{k}^{-1}\mid H^i_{c}({\overline U}, {\Bbb Q}_l))\\
&=\sum_{i=0}^2(-1)^{i}\tr(\varphi_{k}\mid H^i_{\et}({\overline U}, {\Bbb Q}_l(1)))\\
&=\sharp(k)\sum_{i=0}^2(-1)^{i}\tr(\varphi_{k}\mid H^i({\overline U}, {\Bbb Q}_l))\\
&=\sharp(k)\sum_{i=0}^2(-1)^{i}\tr(\varphi_{k}\mid H^i(\Delta_U, {\Bbb Q}_l)).
\endalign
$$
Here, for a profinite group $\Gamma$, we define 
$$H^i(\Gamma, {\Bbb Q}_l)\defeq (\varprojlim H^i(\Gamma, {\Bbb Z}/l^n{\Bbb Z}))
\otimes_{{\Bbb Z}_l}{\Bbb Q}_l,$$
as usual. Thus, $\sharp(X(k))$ and $\sharp(U(k))$ 
can be recovered $\varphi$-group-theoretically. Finally, 
$\sharp(S(k))$ can be recovered as $\sharp(X(k))-\sharp(U(k))$. 

\noindent
(viii) 
First, in the notation of the proof of (iii) above, the image of 
$(\rho^{\det})^{2}=\chi^{2(g+n+b-1)}$ is an open subgroup of the subgroup 
$\overline{\langle p \rangle}$ 
of $({\hat {\Bbb Z}}^{\Sigma^\dag})^{\times}$ (topologically) generated 
by $p$. This characterizes group-theoretically the prime number 
$p$ (in $(\Sigma^\dag)'=\Sigma'\cup\{p\}$), 
by a theorem of Chevalley ([Chevalley], Th\'eor\`eme 1). 
Next, define $m$ to be the minimal positive integer 
with $p^m\in (\rho^{\det})^2(G_k) \  
(\subset ({\hat {\Bbb Z}}^{\Sigma^\dag})^{\times})$. 
Then $\varphi_k$ can be characterized by 
$(\rho^{\det})^2(\varphi_k)=p^m$. 
(See also [Mochizuki2], Remark 1.18.1.) 
\qed\enddemo

\definition{Definition 1.16} 
(i) For each closed subgroup $G$ of $G_{k}$, we denote by 
$k_G$ the subextension of $k$ in $\overline{k}$ 
corresponding to $G$. Observe that, if $G$ is open, 
then $k_G$ is a finite field. 

\noindent
(ii) For each closed subgroup $H$ 
of $\Pi_U$, we set $G_H\defeq \pr_U(H)$ 
and $k_H\defeq k_{G_H}$. 
We denote by 
$U_H$ the (pro-finite, pro-tame, geometrically 
pro-$\Sigma$) covering of $U$ corresponding to $H$. 
Observe that, if $H$ is open, then $U_H$ is 
a hyperbolic curve over the finite field $k_{H}$ 
and $H$ can be identified with $\Pi_{U_H}$. 

\noindent
(iii) Let $H$ be a closed subgroup of $\Pi_U$ and $G$ a 
closed subgroup of $G_H$. Then we set 
$H_G\defeq H\cap\pr_U^{-1}(G)$. Observe 
that $U_{H_G}$ can be identified with $U_H\times_{k_H}k_G$. 

\noindent
(iv) For each open subgroup $H$ of $\Pi_U$, we set 
$$\nu_U(H)\defeq \sharp(U_H(k_{H})).$$
\enddefinition

\proclaim{Corollary 1.17} 
The map $\OSub(\Pi_U) \to {\Bbb Z}_{\geq 0}$, 
$H\mapsto \nu_U(H)$ can be recovered 
$\varphi$-group-theoretically from $\Pi_U$. 
\endproclaim

\demo{Proof}
Since $H=\Pi_{U_H}$, this is immediate from Proposition 1.15(vii). 
\qed
\enddemo

Finally, we shall prove that the set of decomposition groups in $\Pi_U$ 
can be recovered group-theoretically from $\Pi_U$. First, we shall treat 
decomposition groups at points of ${\tilde S}_U$. 

\proclaim{Theorem 1.18} 
{\rm (i)} The set of inertia groups at points of ${\tilde S}_U$ 
(i.e., the image of the injective 
map $I|_{{\tilde S}_U}: {\tilde S}_U \to \Sub(\Delta_U) 
\subset \Sub(\Pi_U)$) can be recovered 
$\varphi$-group-theoretically from $\Pi_U$. 

\noindent
{\rm (ii)} The set of decomposition groups at points of ${\tilde S}_U$ 
(i.e., the image of the injective 
map $D|_{{\tilde S}_U}: {\tilde S}_U \to \Sub(\Pi_U)$) 
can be recovered $\varphi$-group-theoretically from $\Pi_U$. 
\endproclaim

\demo{Proof}
(i) This is due to Nakamura. See [Nakamura1], \S 3 and 
[Nakamura3], 2.1. (See also [Tamagawa1], \S 7, C.)
Strictly speaking, 
Nakamura only treats the case over number fields, but his proof relies 
on Frobenius weights and the same proof goes well over finite fields.

\noindent
(ii) This follows from (i), together with Corollary 1.12(ii). 
\qed
\enddemo

Next, we shall treat decomposition groups at points of ${\tilde U}^{\cl}$. 
This is done along the lines of [Tamagawa1], \S 2, but slightly more subtle 
than the case of [Tamagawa1], due to the existence of the exceptional 
set $E_{\tilde U}$. 

\definition{Definition 1.19} 
(i) We denote by ${\Cal S}(\Pi_U)\  (\subset \Sub(\Pi_U))$ 
the set of closed subgroups $Z$ of $\Pi_U$ 
such that $G_Z$ is open in $G_{k}$ and that 
$\pr_U|_Z:Z\to G_Z$ is an isomorphism. 

\noindent
(ii) For each open subgroup $G$ of $G_{k}$, 
we set 
$${\Cal S}(\Pi_U)_G\defeq\{Z\in {\Cal S}(\Pi_U)\mid G_Z=G\}.$$
Namely, 
${\Cal S}(\Pi_U)_G$ can be identified with the set of group-theoretic 
sections of the surjection $\pr_U|_{(\Pi_U)_G}: (\Pi_U)_G\to G$. 
\enddefinition

\definition{Definition 1.20} 
Let $Z$ be an element of ${\Cal S}(\Pi_U)$. 

\noindent
(i) We define ${\Cal U}(Z)$ to be the set of open subgroups of 
$(\Pi_U)_{G_Z}$ that contain $Z$. 

\noindent
(ii) For each $H\in {\Cal U}(Z)$, we define $m(H,Z)$ to be the number 
of elements $s$ of (a complete system of representatives of) 
$(\Pi_U)_{G_Z}/H$ such that $s^{-1}Zs\subset H$. 
Note that this is a group-theoretic invariant. 

\noindent
(iii) We denote 
by $\nu_\infty(Z)$ the cardinality of $U_Z(k_Z)$. 
(Note that $U_Z(k_Z)$ 
can be identified with the project limit 
of $\{U_H(k_Z)\}_{H\in{\Cal U}(Z)}$.)

\noindent
(iv) We denote by $U_Z(k_{Z})^\ast$ 
the set of points $x$ of $U_Z(k_Z)$ such that the residue field of 
the image of $x$ in $U$ coincides with $k_Z$. (Observe that this 
residue field is included in $k_Z$ in general.) 
We denote by $\nu_\infty^\ast(Z)$ the cardinality of $U_Z(k_Z)^\ast$. 

\noindent
(v) We define 
$U_Z(\overline{k})^{\fin}$ to be the union 
of $U_{Z_{G}}(k_{G})=U_Z(k_{G})$ for all 
open subgroups $G$ of $G_Z$. (N.B. Since $U_Z$ is 
not of finite type over $k_Z$, we have 
$U_Z(\overline{k})^{\fin}\subsetneq U_Z(\overline{k})$.)
We denote by ${\overline\nu}_\infty(Z)$ the cardinality of 
$U_Z(\overline{k})^{\fin}$. Moreover, 
we define $(U_Z)^{\cl,\fin}$ to be the image of 
$U_Z(\overline{k})^{\fin}$ in $(U_Z)^{\cl}$. 
\enddefinition 

\proclaim{Proposition 1.21}
Let $Z$ be an element of ${\Cal S}(\Pi_U)$. 

\noindent
{\rm (i)} Let $G$ be an open subgroup of $G_Z$ and 
$x$ a point of $U_{Z_{G}}(k_G)\subset (U_{Z_{G}})^{\cl}$. Then 
there exists a unique point ${\tilde x}\in{\tilde U}^{\cl}$ above $x$. 
Moreover, $D_{\tilde x}$ contains $Z_{G}$, and, 
in particular, $D_{\tilde x}$ is commensurate to $Z$. 

\noindent
{\rm (ii)} Let ${\tilde x}\in {\tilde U}^{\cl}$ and $x$ the image 
of $\tilde x$ in $(U_Z)^{\cl}$. Then we have
$$x\in U_Z(k_Z) \iff Z\subset D_{\tilde x},$$
$$x\in U_Z(k_Z)^\ast \iff Z= D_{\tilde x},$$
and 
$$x\in (U_Z)^{\cl,\fin}\iff \text{$Z$ and $D_{\tilde x}$ are commensurate.}$$ 

\noindent
{\rm (iii)} We have 
$\nu_\infty^\ast(Z)\leq \nu_\infty(Z)\leq{\overline\nu}_\infty(Z)$ 
and $\nu_\infty(Z)\leq \nu_U((\Pi_U)_{G_Z})<\infty$. 

\noindent
{\rm (iv)} Assume, moreover, that $\Sigma$ is cofinite. Then we have 
${\overline\nu}_\infty(Z)<\infty$. 
\endproclaim

\demo{Proof} 
(i) Take any point ${\tilde x}\in{\tilde U}^{\cl}$ above $x$. 
First note that $D_{\tilde x}\cap Z_{G}$ is the decomposition group 
at ${\tilde x}$ in $Z_{G}$. Thus, since $x$ is $k_{G}$-rational, 
the image of $D_{\tilde x}\cap Z_{G}$ by $\pr_U$ must coincide with $G$. 
Since $\pr_U$ induces an isomorphism $Z_{G}\overset{\sim}\to{\to} G$, 
this implies that $D_{\tilde x}\cap Z_{G}$ coincides with $Z_{G}$. 
It follows from this that $D_{\tilde x}$ contains $Z_{G}$ and that 
there exists only one point (i.e., $\tilde x$) of $\tilde U$ above $x$. 
Finally, since 
$Z_{G}$ is open both in $Z$ and in $D_{\tilde x}$, 
$D_{\tilde x}$ is commensurate to $Z$. (For the latter openness, 
observe that 
$\pr_U$ induces an isomorphism $D_{\tilde x}\overset{\sim}\to{\to} 
\pr_U(D_{\tilde x})$ and that $\pr_U(Z_{G})=G$ is open in 
$\pr_U(D_{\tilde x})$.) 

\noindent
(ii) First, suppose $x\in U_Z(k_Z)$. Then, by (i), $Z\subset D_{\tilde x}$. 
Conversely, suppose $Z\subset D_{\tilde x}$, 
Then the decomposition group $D_{\tilde x} \cap Z$ at ${\tilde x}$ in $Z$ 
coincide with $Z$, which implies $x\in U_Z(k_Z)$. 

Next, we define $x_{U}$ to be the image of $x$ in $U^{\cl}$. 
Suppose $x\in U_Z(k_Z)^\ast$. Then, by (i), 
$Z\subset D_{\tilde x}$. By the definition of 
$U_Z(k_Z)^\ast$, we must have $k(x_U)=k_Z$, or, equivalently, 
$\pr_U(D_{\tilde x})=\pr_U(Z)$. This implies $D_{\tilde x}=Z$. 
Conversely, 
suppose $Z= D_{\tilde x}$. Then 
$\pr_U(D_{\tilde x})=\pr_U(Z)$, or, equivalently, $k(x_U)=k_Z$. This 
implies $x\in U_Z(k_Z)^\ast$. 

Finally, for each open subgroup $G$ of $G_Z$, denote by $x_G$ 
the image of $\tilde x$ in $(U_{Z_G})^{\cl}$. Then 
$$
\align
&x\in (U_Z)^{\cl,\fin}\\
\iff& \text{$x_G\in U_{Z_G}(k_{Z_G})$ for some open subgroup $G$ of $G_Z$}\\
\iff& \text{$Z_G\subset D_{\tilde x}$ for some open subgroup $G$ of $G_Z$}\\
\iff&\text{$Z$ and $D_{\tilde x}$ are commensurate.}
\endalign$$ 

\noindent
(iii) The first two inequalities are clear. To see the third inequality, 
it suffices to prove that the natural 
map $U_{Z}(k_Z)\to U_{(\Pi_U)_{G_Z}}(k_Z)$ is injective. 
For this, take $x,x'\in U_{Z}(k_Z)$ and suppose that the images of 
$x,x'$ in $U_{(\Pi_U)_{G_Z}}(k_Z)$ coincide with each other. 
Take the unique points ${\tilde x},{\tilde x}'\in {\tilde U}^{\cl}$ 
above $x,x'$, respectively. Then, by (i), $D_{\tilde x}$ 
and $D_{{\tilde x}'}$ are commensurate to each other. 
On the other hand, since the images of ${\tilde x}, {\tilde x}'$ in 
$(U_{(\Pi_U)_{G_Z}})^{\cl}$ coincide with each other and are $k_Z$-rational, 
we see that their images in ${\overline U}^{\cl}$ must coincide 
with each other. It follows from these and Proposition 1.8(iii) that 
${\tilde x}={\tilde x}'$ holds. Thus, $x=x'$ holds, as desired. 
Finally, the fourth equality is clear. 

\noindent
(iv) This follows from (ii) and Proposition 1.8(vii). 
(Observe that the natural surjective map 
$U_Z(\overline{k})^{\fin}\to(U_Z)^{\cl,\fin}$ 
is quasi-finite.) \qed
\enddemo

\proclaim{Corollary 1.22} 
Let $Z$ be an element of ${\Cal S}(\Pi_U)$. Then we have 

\noindent
{\rm (i)} 
There exists ${\tilde x}\in {\tilde U}^{\cl}$ such that $Z=D_{\tilde x}$
(respectively, $Z\subset D_{\tilde x}$, respectively, $Z$ is commensurate to $D_{\tilde x}$), 
if and only if $\nu_\infty^\ast(Z)>0$ 
(respectively, $\nu_{\infty}(Z)>0$, respectively, ${\overline\nu}_\infty(Z)>0$). 

\noindent
{\rm (ii)} 
There exist more than one ${\tilde x}\in {\tilde U}^{\cl}$ such that 
$Z$ is commensurate to $D_{\tilde x}$, 
if and only if ${\overline\nu}_\infty(Z)>1$. 

\noindent
{\rm (iii)} There exists ${\tilde x}\in {\tilde U}^{\cl}-E_{\tilde U}$ 
(respectively, ${\tilde x}\in E_{\tilde U}$) 
such that $Z=D_{\tilde x}$
if and only if ${\overline\nu}_{\infty}(Z)=\nu_\infty^\ast(Z)=1$ 
(respectively, $\nu^\ast_{\infty}(Z)>0$ and ${\overline\nu}_\infty(Z)>1$). 
\endproclaim

\demo{Proof}
(i) This is immediate from Proposition 1.21(ii). 

\noindent
(ii) By definition, ${\overline\nu}_\infty(Z)>1$ if and only if 
$\nu_{\infty}(Z_G)>1$ for some open subgroup $G$ of $G_Z$. Thus, 
the assertion follows from (the first statement of) 
Proposition 1.21(ii) and (the uniqueness statement of) 
Proposition 1.21(i). 

\noindent
(iii) It follows formally from (i) and (ii) 
that $\nu^\ast_{\infty}(Z)>0$ and ${\overline\nu}_\infty(Z)>1$ 
(respectively, $\leq 1$)
if and only if $Z=D_{\tilde x}$ for some 
${\tilde x}\in{\tilde U}^{\cl}$ and 
$Z$ is commensurate to $D_{\tilde x}$ for more than (respectively, at most) 
one ${\tilde x}\in{\tilde U}^{\cl}$. This last statement 
is equivalent to saying that $Z=D_{\tilde x}$ for some 
${\tilde x}\in E_{\tilde U}$ (respectively, ${\tilde x}\in {\tilde U}^{\cl}- 
E_{\tilde U}$). 
This, together with Proposition 1.21(iii) 
(or, more specifically, 
the fact $\nu_{\infty}^\ast(Z)\leq{\overline\nu}_{\infty}(Z)$), 
completes the proof. 
\qed
\enddemo

\proclaim{Proposition 1.23} 
Let $Z$ be an element of ${\Cal S}(\Pi_U)$. 

\noindent
{\rm (i)} We have 
$$\lim\Sb H\in{\Cal U}(Z) \\ H\to Z\endSb \frac{\nu_U(H)}{m(H,Z)}
=\nu_{\infty}(Z).$$
More precisely, there exists $H_0\in{\Cal U}(Z)$ such that, for any 
$H\in{\Cal U}(Z)$ with $H\subset H_0$, we have 
$$\frac{\nu_U(H)}{m(H,Z)}=\nu_{\infty}(Z).$$

In particular, $\nu_{\infty}(Z)$ is a $\varphi$-group-theoretic invariant. 

\noindent
{\rm (ii)} 
Set $C\defeq C_{\Pi_U}(Z)$, which is isomorphic to 
${\hat {\Bbb Z}}$. Then we have 
$$\nu_{\infty}^\ast(Z)=\sum_{d\mid N} \mu(N/d) \nu_{\infty}(C^d),$$
where $N\defeq (C:Z)$, 
$C^d\defeq\{\sigma^d\mid \sigma\in C\}$, 
and $\mu$ stands for M\"obius' function. 

In particular, $\nu_{\infty}^\ast(Z)$ is a 
$\varphi$-group-theoretic invariant. 

\noindent
{\rm (iii)} We have 
$$\sup\Sb G\in\OSub(G_Z)\endSb 
\nu_{\infty}(Z_{G})={\overline\nu}_{\infty}(Z).$$

In particular, ${\overline\nu}_{\infty}(Z)$ is a 
$\varphi$-group-theoretic invariant. 
\endproclaim

\demo{Proof}
(i) We define $U_{(\Pi_U)_{G_Z}}(k_Z)^{\infty}$ 
to be the image of 
$U_Z(k_Z)$ in $U_{(\Pi_U)_{G_Z}}(k_Z) (= U(k_Z))$. 
On one hand, the proof of (the third inequality of) 
Proposition 1.21(iii) shows that 
the natural surjection $U_{Z}(k_Z)\to U_{(\Pi_U)_{G_Z}}(k_Z)^{\infty}$ 
is a bijection. On the other hand, 
since 
$U_Z(k_Z)=\varprojlim_{H\in{\Cal U}(Z)}U_H(k_Z)$ 
and $\sharp(U_H(k_Z))<\infty$ for each $H\in{\Cal U}(Z)$, we see that 
there exists $H_0\in{\Cal U}(Z)$ such that 
$U_{(\Pi_U)_{G_Z}}(k_Z)^{\infty}$ coincides with 
the image of $U_{H_0}(k_Z)$ in $U_{(\Pi_U)_{G_Z}}(k_Z)$. 

Take any $H\in{\Cal U}(Z)$ with $H\subset H_0$. Then 
each point of $U_H(k_Z)$ is above some point of 
$U_{(\Pi_U)_{G_Z}}(k_Z)^{\infty}$. Moreover, for each point 
$x\in U_Z(k_Z)$, take the unique point ${\tilde x}\in {\tilde U}^{\cl}$ 
above $x$. Then, by Proposition 1.21(i), the decomposition group at 
$\tilde x$ in $(\Pi_U)_{G_Z}$ coincides with $Z$. From this, we see that 
$m(H,Z)$ is defined so as to coincide with the cardinality of 
the fiber of the map $U_H(k_Z)\to U_{(\Pi_U)_{G_Z}}(k_Z)$ at $x_{(\Pi_U)_{G_Z}}$, 
where $x_{(\Pi_U)_{G_Z}}$ is the image of 
$x$ in $U_{(\Pi_U)_{G_Z}}(k_Z)$. From these, we conclude that the quantity 
$\nu_U(H)/m(H,Z)$ coincides with the cardinality 
$\sharp(U_{(\Pi_U)_{G_Z}}(k_Z)^{\infty})$, 
as desired. 

\noindent
(ii) 
First, by Lemma 1.7, we see that $C\in{\Cal S}(\Pi_U)$ and
that $C$ is isomorphic to $\hat{\Bbb Z}$. 

Let $x$ be a point of $U_Z(k_Z)$. We claim that 
$x\not\in U_Z(k_Z)^\ast$ if and only if there exists 
$Z'\in{\Cal S}(\Pi_U)$ with $Z'\supsetneq Z$, such that 
the image in $(U_{Z'})^{\cl}$ of $x\in U_Z(k_Z)\subset (U_Z)^{\cl}$ 
is $k_{Z'}$-rational. Indeed, to see the `if' part, 
observe that the natural morphism $U_Z\to U$ factors through 
$U_Z \to U_{Z'}$. Thus, if the image of $x$ in 
$(U_{Z'})^{\cl}$ is $k_{Z'}$-rational, so is the image of $x$ in 
$U^{\cl}$, hence $x\not\in U_Z(k_Z)^\ast$. Conversely, suppose 
$x\not\in U_Z(k_Z)^\ast$ and take 
the unique point ${\tilde x}\in {\tilde U}^{\cl}$ above $x$. 
As the residue field of the image of $x$ in $U$ is strictly smaller 
than $k_Z$, the image of $D_{\tilde x}$ in $G_{k}$ must be strictly 
larger than $G_Z$. Now, it is easy to see from this that 
$Z'\defeq D_{\tilde x}$ has the desired property. 

Now, consider $Z'\in{\Cal S}(\Pi_U)$ with $Z'\supset Z$. 
Then we have $Z\subset Z'\subset C$, 
which implies that $Z'=C^d$ for some (unique) 
$d$ dividing $N$. 
We see $U_{Z'}= U_C \times_{k_{C}}k_{Z'}$, and, 
in particular, $U_{Z}= U_C \times_{k_{C}}k_{Z}$. 
Thus, 
the image of $x$ in $(U_{Z'})^{\cl}$ is $k_{Z'}$-rational 
if and only if $x\in U_Z(k_Z)=U_{C}(k_Z)$ falls in 
$U_{Z'}(k_{Z'})=U_C(k_{Z'})$. 

These, together with the so-called inclusion-exclusion principle, 
imply the desired formula. 

\noindent
(iii) Immediate from the definitions. 
\qed
\enddemo

\proclaim{Theorem 1.24} 
The set of decomposition groups at points of ${\tilde U}^{\cl}$ 
(respectively, ${\tilde U}^{\cl}- E_{\tilde U}$, respectively, $E_{\tilde U}$) 
(i.e., the image of the 
map 
$D|_{{\tilde U}^{\cl}}: {\tilde U}^{\cl} \to \Sub(\Pi_U)$ 
(respectively, 
$D|_{{\tilde U}^{\cl}-E_{\tilde U}}: {\tilde U}^{\cl}-E_{\tilde U} \to \Sub(\Pi_U)$, 
respectively, 
$D|_{E_{\tilde U}}: E_{\tilde U} \to \Sub(\Pi_U)$))
can be recovered $\varphi$-group-theoretically from $\Pi_U$. 
\endproclaim

\demo{Proof} This follows formally from Corollary 1.22 and Proposition 1.23. 
\qed\enddemo

\proclaim{Corollary 1.25} 
The set of decomposition groups at points of $({\tilde X}_U)^{\cl}$ 
(i.e., the image of the 
map 
$D: ({\tilde X}_U)^{\cl} \to \Sub(\Pi_U)$) 
can be recovered $\varphi$-group-theoretically from $\Pi_U$. 
\endproclaim

\demo{Proof} This is immediate from Theorem 1.18(ii) and Theorem 1.24. 
(See also [Mochizuki2], Remark 1.18.2.) 
\qed\enddemo

\subhead
\S 2. 
Cuspidalizations of Proper Hyperbolic Curves
\endsubhead

In this \S, we review the main results of Mochizuki's theory 
of cuspidalizations of fundamental groups of proper hyperbolic curves, 
developed in [Mochizuki2], which plays an important role in this paper. 
We maintain the notations of \S 1 and further assume 
$X=U$. (Thus, the finite set $S$ in \S 1 is empty, and, in this \S, 
we save the symbol $S$ for another finite set of closed points of $X$.) 
Accordingly, $X$ is a proper hyperbolic curve over a finite field 
$k=k_X$. 



Recall that $\Delta _X$ stands for the maximal 
pro-$\Sigma$ quotient of $\pi _1(\overline X)$, that 
$\Pi _X$ stands for 
$\pi _1(X)/\Ker (\pi _1(\overline X)\twoheadrightarrow \Delta_X)$, and that 
they fit into the following exact sequence: 
$$1\to \Delta _X \to \Pi _X@>\pr_X>> G_k\to 1.$$

Similarly, if we write $X\times X\overset \text {def}\to=X\times _kX$, 
then we obtain (by considering the maximal pro-$\Sigma$ quotient $\Delta _{X\times X}$ 
of $\pi _1(\overline {X\times X}$)) an exact sequence:
$$1\to \Delta _{X\times X} \to \Pi _{X\times X}\to G_k\to 1,$$
where $\Pi _{X\times X}$ (respectively, $\Delta _{X\times X}$) may be 
identified with $\Pi _X\times _{G_k}\Pi_X$ (respectively, $\Delta _X\times 
\Delta_X$).



\definition{Definition 2.1} (cf. [Mochizuki2], Definition 1.5(i).) 
Let $H$ be a profinite 
group equipped with a homomorphism $H\to \Pi _X$.
Then we shall refer to the kernel $I_H$ of $H\to 
\Pi _X$ as the cuspidal subgroup of $H$ (relative to $H\to \Pi _X$).
We shall refer to an inner automorphism of $H$ by an element of $I_H$ as a
cuspidally inner automorphism.
We shall say that $H$ is cuspidally abelian (respectively, 
cuspidally pro-$\Sigma ^{*}$, where $\Sigma ^{*}$ is a set of prime numbers)
(relative to $H\to \Pi _X$) if $I_H$ is abelian (respectively, if $I_H$ is a 
pro-$\Sigma ^{*}$ group). If $H$ is cuspidally abelian, then observe that 
$H/I_H$ acts naturally (by conjugation) on $I_H$. We shall say that $H$ is 
cuspidally central (relative to $H\to \Pi _X$) if this action of 
$H/I_H$ on $I_H$ is trivial. Also, we shall use the same terminology 
for  $H\to \Pi _X$ when $\Pi _X$ is replaced by $\Delta _X$, 
$\Pi _{X\times X}$, or $\Delta _{X\times X}$.
\enddefinition

For a finite subset $S\subset X^{\cl}$ write $U_S\overset \text {def}
\to=X-S$. Let $\Delta _{U_S}$ be the maximal cuspidally
(relative to the natural map to $\Delta _X$) 
pro-$\Sigma ^{\dag}$ quotient of the maximal pro-$\Sigma$ quotient of 
the tame fundamental group of $\overline {U_{S}}$ (where ``tame'' is 
with respect to the complement of $U_S$ in $X$), and let $\Pi 
_{U_S}$ be the corresponding quotient $\pi _1(U_S)/\Ker (\pi _1
(\overline {U_S})\twoheadrightarrow \Delta _{U_S}$) of $\pi_1(U_S)$. 
Thus, we have an exact sequence:  
$$1\to \Delta _{U_S}\to \Pi _{U_S} \to G_k\to 1,$$
which fits into the following commutative diagram: 
$$\matrix
1&\to &\Delta _{U_S}&\to &\Pi _{U_S} &\to &G_k&\to& 1\\
&&&&&&&&\\
&& \downarrow && \downarrow && \Vert && \\
&&&&&&&&\\
1&\to &\Delta _{X}&\to &\Pi _{X} &\to& G_k&\to& 1.
\endmatrix
$$

Further, let $\iota :X \to X\times X$ be the diagonal
morphism, and write
$$U_{X\times X}\overset \text {def}\to=X\times X-\iota (X).$$
We shall denote by $\Delta _{U_{X\times X}}$ the maximal cuspidally 
(relative to the natural map to $\Delta _{X\times X}$) 
pro-$\Sigma^{\dag}$ quotient of the maximal pro-$\Sigma$ quotient of 
the tame fundamental group of 
$(U_{X\times X})_{\bar k}$ (where ``tame'' is 
with respect to the divisor $\iota (X)\subset X\times X$), and by $\Pi 
_{U_{X\times X}}$ the corresponding quotient $\pi _1(U_{X\times X})/\Ker (\pi _1
(\overline {U_{X\times X}})\twoheadrightarrow \Delta _{U_{X\times X}}$) of 
$\pi _1(U_{X\times X})$.
Thus, we have an exact sequence:
$$1\to \Delta _{U_{X\times X}}\to \Pi _{U_{X\times X}}\to G_k\to 1,$$
which fits into 
the following commutative diagram: 
$$\matrix
1&\to &\Delta _{U_{X\times X}}&\to &\Pi _{U_{X\times X}}&\to& G_k&\to& 1\\
&&&&&&&&\\
&&\downarrow&&\downarrow&&\Vert&&\\
&&&&&&&&\\
1&\to &\Delta _{X\times X} &\to &\Pi _{X\times X}&\to &G_k&\to& 1. 
\endmatrix
$$

Finally, set
$$M_X\overset \text {def}\to=\Hom _ {\hat \Bbb Z^{\Sigma^{\dag}}}(H^2(\Delta_X,
\hat \Bbb Z^{\Sigma^{\dag}}),\hat \Bbb Z^{\Sigma^{\dag}}).$$
Thus, $M_X$ is a free $\hat \Bbb Z^{\Sigma^{\dag}}$-module of rank $1$, 
and $M_X$
is isomorphic to $\hat \Bbb Z^{\Sigma^{\dag}}(1)$ as a $G_k$-module 
(where the ``$(1)$'' denotes a ``Tate twist'', 
i.e., $G_k$ acts on $\hat \Bbb Z^{\Sigma^{\dag}}
(1)$ via the cyclotomic character) (cf. [Mochizuki2], the discussion 
following Proposition 1.1). 

For the rest of this \S, let $X$, $Y$ be proper, 
hyperbolic curves over finite fields $k_X$, $k_Y$ of characteristic
$p_X$, $p_Y$, respectively. Let $\Sigma_X$ (respectively, $\Sigma_Y$)
be a set of prime numbers that contains at least one prime number
different from $p_X$ (respectively, $p_Y$).
Write $\Delta _X$ (respectively, $\Delta _Y$) 
for the maximal pro-$\Sigma _X$ quotient of $\pi _1(\overline X)$ 
(respectively, the maximal  pro-$\Sigma_Y$ quotient of $\pi _1
(\overline Y)$), and $\Pi _X$ (respectively, $\Pi _Y$) for the 
quotient $\pi_1(X)/\Ker (\pi_1(\overline X)\twoheadrightarrow \Delta_X)$
of $\pi _1(X)$ (respectively, the quotient 
$\pi_1(Y)/\Ker (\pi_1(\overline Y)\twoheadrightarrow \Delta_Y)$ of $\pi _1(Y)$). 

Let 
$$\alpha :\Pi _X\isom \Pi _Y$$
be an isomorphism of profinite groups. 

%
%
%

The following is one of the main results of Mochizuki's theory 
(cf. [Mochizuki2], Theorem 1.16(iii)).

\proclaim{Theorem 2.2 (Reconstruction of Maximal Cuspidally Abelian Extensions)} 
Let $\iota _X:X\to X\times X$ (respectively, 
$\iota _Y:Y\to Y\times Y$) be the diagonal morphism, and write
$U_{X\times X}\overset \text {def}
\to=X\times X- \iota (X)$ (respectively, $U_{Y\times Y}
\overset \text {def}\to=Y\times Y- \iota (Y)$). Denote by 
$\Pi _{U_{X\times X}}\twoheadrightarrow\Pi _{U_{X\times X}}^{\text{\rm c-ab}}$, 
$\Pi _{U_{Y\times Y}}\twoheadrightarrow\Pi _{U_{Y\times Y}}^{\text{\rm c-ab}}$ 
the maximal cuspidally (relative to 
the natural surjections  $\Pi _{U_{X\times X}}\twoheadrightarrow \Pi_
{X\times X}$, $\Pi _{U_{Y\times Y}}\twoheadrightarrow \Pi_{Y\times Y}$, 
respectively)
abelian quotients. Then there is a commutative diagram:
$$
\CD
\Pi _{U_{X\times X}}^{\text{\rm c-ab}} @>{\alpha^{\text{\rm c-ab}}}>> \Pi _{U_{Y\times Y}}^{\text{\rm 
c-ab}} \\
@VVV                @VVV \\
\Pi_{X\times X}   @>{\alpha\times \alpha}>>  \Pi _{Y\times Y}\\
\endCD
$$
where $\alpha^{\text{\rm c-ab}}$ is an isomorphism which is well-defined up to 
cuspidally inner automorphism (i.e., an inner automorphism of $\Pi _ 
{U_{Y\times Y}}^{\text{\rm c-ab}}$ by an element of the cuspidal subgroup 
$\Ker (\Pi _{U_{Y\times Y}}^{\text{\rm c-ab}}\twoheadrightarrow 
\Pi _{Y\times Y})$). Moreover, the correspondence 
$$\alpha \mapsto \alpha ^{\text{\rm c-ab}}$$
is functorial (up to cuspidally inner automorphism) with respect to $\alpha$.
\endproclaim

\demo {Proof}
See [Mochizuki2], Theorem 1.16(iii). 
\qed
\enddemo

In this and next \S\S, we sometimes refer to the decomposition group
$D_{\tilde x}$ as the decomposition group of $\Pi_X$ at $x$, and denote it simply
by $D_x$. Thus, $D_x$ is well-defined only up to conjugation by an element of $\Pi_X$. 

For the rest of this \S, we shall assume 
that $\alpha$ is Frobenius-preserving (cf. Definition 1.14).
(Note that this assumption is automatically satisfied 
in the case where $\Sigma_X$ 
and $\Sigma _Y$ are cofinite by Theorem 1.15(viii).) 
Thus, by Theorem 1.24, 
one deduces naturally from $\alpha$ a bijection
$$\phi:X^{\cl}-E_X\isom Y^{\cl}-E_Y$$ 
such that 
$$\alpha (D_x)=D_{\phi (x)}$$ 
holds (up to conjugation) for any $x\in X^{\cl}- E_X$. 
(Note that $E_X$ (respectively, $E_Y$) is a finite set by Proposition 1.8(vi),  
if $\Sigma_X$ (respectively, $\Sigma_Y$) is cofinite.) 

As an important consequence of Theorem 2.2 we deduce the following:

\proclaim{Corollary 2.3} With the above assumptions,
let $S\subset X^{\cl}-E_X$,  
$T\subset Y^{\cl}-E_Y$ be finite subsets that correspond 
to each other via $\phi$. Then $\alpha$, $\alpha ^{\text{\rm c-ab}}$ induce 
isomorphisms (well-defined up to cuspidally inner automorphisms, i.e.,
inner automorphisms by elements of $\Ker (\Pi _{V_T}^{\text{\rm c-ab}} \to \Pi_Y)$)
$$\alpha_{S,T}^{\text{\rm c-ab}}: 
\Pi_{U_S}^{\text{\rm c-ab}}\isom \Pi _{V_T}^{\text{\rm c-ab}}$$
lying over $\alpha$, where $U_S\defeq X-S$, $V_T\defeq Y-T$, 
and $\Pi_{U_S}\twoheadrightarrow \Pi_{U_S}^{\text{\rm c-ab}}$,
$\Pi_{V_T}\twoheadrightarrow \Pi_{V_T}^{\text{\rm c-ab}}$, are the maximal cuspidally 
abelian quotients (relative to the maps $\Pi_{U_S}\twoheadrightarrow \Pi_X$,
$\Pi_{V_T}\twoheadrightarrow \Pi_Y$, respectively).
These isomorphisms are functorial with respect to $\alpha$, $S$, 
$T$, as well as with respect to passing to connected finite \'etale coverings
of $X$, $Y$, which arise from open subgroups of $\Pi_X$, $\Pi_Y$, 
in the following sense:
Let $\xi:X'\to X$ (respectively, $\eta : Y'\to Y$) be a finite \'etale covering
which arises from the open subgroup $\Pi_{X'}\subseteq \Pi_X$ (respectively,  
$\Pi_{Y'}\subseteq \Pi_Y$), such that $\alpha (\Pi_{X'})=\Pi_{Y'}$; 
set $U'_{S'}\defeq X'-S'$, $V'_{T'}\defeq Y'-T'$, 
$S' \defeq \xi ^{-1} (S)$, $T' \defeq \eta^{-1} (T)$; 
and denote by $\alpha'$ the isomophism 
$\Pi_{X'}\isom\Pi_{Y'}$ induced by $\alpha$. 
Then we have the following commutative diagram:
$$
\CD
\Pi_{U'_{S'}}^{\text{\rm c-ab}} @>{(\alpha')_{S',T'}^{\text{\rm c-ab}}}>> 
\Pi_{V'_{T'}}^{\text{\rm c-ab}}  \\
@VVV                            @VVV        \\
\Pi_{U_S}^{\text{\rm c-ab}} @>{\alpha_{S,T}^{\text{\rm c-ab}}}>> \Pi_{V_{T}}^{\text{\rm c-ab}}  \\
\endCD
$$
where the vertical arrows are the natural maps. 
\endproclaim

\demo {Proof} 
The proof of [Mochizuki2], Theorem 2.5(i) 
(where $E_X=E_Y=\emptyset$ is assumed) 
works as it is. 
\qed
\enddemo

Next, let 
$$1\to M_X\to \Cal D\to \Pi_{X\times X}\to 1$$
be a fundamental extension, i.e., an extension whose corresponding 
extension class in $H^2_{\et}(X\times X,M_X)$ (via the natural identification
$H^2(\Pi_{X\times X}, M_X)\isom H^2_{\et}(X\times X,M_X)$ 
(cf. [Mochizuki2], Proposition 1.1)) coincides with the first 
(\'etale) Chern class of the diagonal $\iota (X)$ (cf. [Mochizuki2], Proposition 1.6). 
Let $x,y\in X(k)$ and 
write $D_x,D_y\subset \Pi_X$ for the associated 
decomposition groups (which are well-defined up to conjugation). Set
$$\Cal D_x\overset \text {def}\to=\Cal D|D_x\times _{G_k}\Pi_X,\ \ \ 
\Cal D_{x,y}\overset \text {def}\to=\Cal D|D_x\times _{G_k}D_y.$$
Thus, $\Cal D_x$ (respectively, $\Cal D_{x,y}$) is an extension of $\Pi_X$
(respectively, $G_k$) by $M_X$. Similarly, if $D=\sum _{i}m_i.x_i,\  
E=\sum_jn_j.y_j$ are divisors on $X$ supported 
on $k$-rational points, then set 
$$\Cal D_D\overset \text {def}\to=\sum _i m_i.\Cal D_{x_i},\  
\Cal D_{D,E}\overset \text {def}\to=\sum _{i,j} m_i.n_j.\Cal D_{x_i,y_j}$$
where the sums are to be understood as sums of extensions of $\Pi _X$, 
$G_k$, respectively, 
by $M_X$, i.e., the sums are induced by the additive structure of $M_X$.

For a finite subset $S\subset X(k)$, we shall write
$$\Cal D_S\overset \text {def}\to=\underset {x\in S}\to \prod \Cal D_x$$
where the product is to be understood as a fiber product over $\Pi_X$.
Thus, $\Cal D_S$ is an extension of $\Pi_X$ by a product of copies of $M_X$ 
indexed by the points of $S$. We shall refer to  $\Cal D_S$ as the 
$S$-cuspidalization of $\Pi_X$. 
Observe that if $T\subset X(k)$
is a finite subset containing $S$, then we obtain a natural  
projection morphism $\Cal D_T\to\Cal D_S$. More generally, for a 
finite subset $S\subset X^{\cl}$ which does not 
necessarily consist of $k$-rational
points, one can still construct the object $\Cal D_S$ by passing to a finite 
extension $k_S$ of $k$ over which the points of $S$ are rational, performing the
above construction over $k_S$, and then descending to $k$. 
(See [Mochizuki2], Remark 1.10.1 for more details.) 

\proclaim {Proposition 2.4 
(Maximal Geometrically Cuspidally Central Quotients)} 

\noindent
{\rm (i)}\ For $S\subset X^{\cl}$ a finite subset, the $S$-cuspidalization 
$\Cal D_S$ of $\Pi_X$ may be identified with the 
quotient $\Pi_{U_S}\twoheadrightarrow \Pi_{U_S}^{\text{\rm c-cn}}
\overset \text {def}\to=\Pi_{U_S}/ \Ker (\Delta _{U_S}\twoheadrightarrow \Delta _{U_S}^{\text{\rm 
c-cn}})$
of $\Pi_{U_S}$, where $\Delta _{U_S}^{\text{\rm c-cn}}$ is the maximal cuspidally
central quotient of $\Delta _{U_S}$ relative to the natural map 
$\Delta _{U_S} \twoheadrightarrow \Delta _X$.

\noindent
{\rm (ii)}\ The fundamental extension $\Cal D$ may be identified with the 
quotient $\Pi_{X\times X}\twoheadrightarrow 
\Pi_{U_{X\times X}}^{\text{\rm c-cn}}\overset \text {def}\to=\Pi_{U_{X\times X}}/\Ker 
(\Delta _{U_{X\times X}}\twoheadrightarrow   \Delta _{U_{X\times X}}
^{\text{\rm c-cn}})$ of $\Pi_{U_{X\times X}}$, where $\Delta _{U_{X\times X}}^{\text{\rm c-cn}}$ is
the maximal cuspidally central quotient of $\Delta _{U_{X\times X}}$
relative to the natural map $\Delta _{U_{X\times X}}\twoheadrightarrow \Delta_{X\times X}$.
\endproclaim

\demo 
{Proof}\ See [Mochizuki2], Proposition 1.8(iii)(iv). 
(Precisely speaking, Proposition 1.8(iii) loc. cit. only treats the 
special case where $S\subset X(k)$ holds. However, the proof for 
the general case is easily reduced to this special case by passing 
to a finite extension of $k$. cf. Remark 1.10.1 loc. cit.) 
\qed
\enddemo

\definition {Remark 2.5}
Let $\Cal D$ (respectively, $\Cal E$) be a fundamental extension 
of $X$ (respectively, $Y$). 
The isomorphism $\alpha:\Pi_X\isom \Pi_Y$ induces an isomorphism:
$$\Cal D\isom \Cal E$$
up to cyclotomically inner automorphisms (i.e., inner automorphisms 
by elements of $M_X, M_Y$) and the actions of $(k_X^{\times})^{\Sigma_X^{\dag}}, 
(k_Y^{\times})^{\Sigma_Y^{\dag}}$, where $(k_X^{\times})^{\Sigma_X^{\dag}}$ 
(respectively, $(k_Y^{\times})^{\Sigma_Y^{\dag}}$) is the maximal $\Sigma_X^{\dag}$-
(respectively, $\Sigma_Y^{\dag}$-) quotient of $k_X^{\times}$ (respectively, 
$k_Y^{\times}$) (cf. [Mochizuki2], Proposition 1.4(ii)). 
Moreover, let $S\subset X^{\cl}-E_X$ and $T\subset Y^{\cl}-E_Y$ be 
as in Corollary 2.3 and write $\Cal D_S$ (respectively, $\Cal E_T$) 
for the $S$-cuspidalization of $\Pi_X$ 
(respectively, the $T$-cuspidalization of $\Pi_Y$). 
Then the isomorphism $\Cal D\isom \Cal E$ induces an isomophism 
$$\Cal D_S\isom \Cal E_T$$
lying over $\alpha$. 

On the other hand, 
let $\Pi_{U_S}\twoheadrightarrow \Pi_{U_S}^{\text{\rm c-cn}}$ and
$\Pi_{V_T}\twoheadrightarrow \Pi_{V_T}^{\text{\rm c-cn}}$ be 
the maximal geometrically cuspidally central quotients 
(here, $U_S\overset \text {def}\to=X-S$, 
$V_T\overset \text {def}\to=Y-T$) (cf. Proposition 2.4). 
Note that the isomorphism 
$\alpha_{S,T}^{\text{\rm c-ab}}: 
\Pi_{U_S}^{\text{\rm c-ab}}\isom \Pi _{V_T}^{\text{\rm c-ab}}$ 
in Corollary 2.3 naturally induces an isomorphism
$$\Pi_{U_S}^{\text{\rm c-cn}}\isom \Pi _{V_T}^{\text{\rm c-cn}}$$ 
lying over $\alpha$, 
which is well-defined up to cuspidally inner automorphism. 
Now, by Proposition 2.4(i), $\Pi_{U_S}^{\text{\rm c-cn}}$ 
(respectively, $\Pi_{V_T}^{\text{\rm c-cn}}$) may be identified with 
$\Cal D_S$ (respectively, $\Cal E_T$). 
Thus, we deduce another isomorphism 
$$\Cal D_S\isom \Cal E_T$$
lying over $\alpha$. 

Now, the above two isomorphisms between $\Cal D_S$ and $\Cal E_T$ 
coincide with each other up to cyclotomically inner automorphisms 
and the actions of $(k_X^{\times})^{\Sigma_X^{\dag}}, 
(k_Y^{\times})^{\Sigma_Y^{\dag}}$.
\enddefinition

Another main result of Mochizuki's theory is the following, 
which allows us to recover $\varphi$-group-theoretically the maximal cuspidally 
pro-$l$ extension of $\Pi_X$, in the case where 
the set of cusps consists of a single rational point. 

\proclaim {Theorem 2.7 (Reconstruction of Maximal Cuspidally Pro-$l$ Extensions)} 

\noindent
Let $x_*\in X(k_X)$, $y_*\in Y(k_Y)$, and 
set $S\overset \text {def}\to=\{x_*\}$, $T\overset \text {def}\to=\{y_*\}$,
$U_S\overset \text {def}\to=X-S$, $V_T\overset \text {def}\to=Y-T$. 
Assume that the Frobenius-preserving isomorphism $\alpha:\Pi_X\isom\Pi_Y$ 
maps 
the decomposition group of $x_*$ in $\Pi_X$ (which is well-defined up to conjugation) 
to the decomposition group of $y_*$ in $\Pi_Y$ (which is well-defined up to 
conjugation). Then, for 
each prime $l\in \Sigma ^{\dag}$ (thus, $l\neq p$), 
there exists a commutative diagram:
$$
\CD
\Pi_{U_S}^{\text{\rm c-$l$}} @>\alpha^{\text{\rm c-$l$}}>>  
\Pi_{V_T}^{\text{\rm c-$l$}}\\
@VVV                   @VVV   \\
\Pi_X        @>\alpha>> \Pi_Y \\
\endCD
$$ 
in which $\Pi_{U_S}\twoheadrightarrow  \Pi_{U_S}^{\text{\rm c-$l$}}$, $\Pi_{V_T}
\twoheadrightarrow  \Pi_{V_T}^{\text{\rm c-$l$}}$ are the maximal cuspidally pro-$l$ 
quotients (relative to the maps $\Pi_{U_S}\twoheadrightarrow    \Pi_X$, 
$\Pi_{V_T}\twoheadrightarrow  \Pi_Y$, respectively), the vertical arrows are the 
natural surjections, and $\alpha^{\text{\rm c-$l$}}$
is an isomorphism well-defined up to composition with a cuspidally inner 
automorphism (i.e., an inner automorphism by an element of 
$\Ker (\Pi_{V_T}^{\text{\rm c-$l$}} \to \Pi_Y$)), which is compatible relative to the natural 
surjections
$$\Pi_{U_S}^{\text{\rm c-$l$}}\twoheadrightarrow \Pi_{U_S}^{\text{\rm c-ab},l}
,\ \ \ \ \Pi_{V_T}^{\text{\rm c-$l$}}\twoheadrightarrow \Pi_{V_T}^{\text{\rm c-ab},l}$$
where the subscript ``$\text{\rm c-ab},l$'' denotes the maximal cuspidally pro-$l$ 
abelian quotient, with the isomorphism
$$\alpha_{S,T}^{\text{\rm c-ab}}: 
\Pi_{U_S}^{\text{\rm c-ab}}\isom \Pi _{V_T}^{\text{\rm c-ab}}$$
in Corollary 2.3. Moreover, $\alpha^{\text{\rm c-$l$}}$ is compatible, up to cuspidally
inner automorphisms, with the decomposition groups of $x_*$, $y_*$ in 
$\Pi_{U_S}^{\text{\rm c-$l$}}$, $\Pi_{V_T}^{\text{\rm c-$l$}}$.
\qed
\endproclaim

\demo
{Proof}\ See [Mochizuki2], Theorem 3.10. 
\qed
\enddemo

\subhead
\S 3. Kummer Theory and Anabelian Geometry
\endsubhead

We maintain the notations of \S 2. 
If $n$ is 
an integer all of whose prime factors belong to $\Sigma ^{\dag}$, 
then we have the Kummer exact sequence
$$1\to \mu_n\to \Bbb G_m\to \Bbb G_m\to 1,$$
where $\Bbb G_m\to \Bbb G_m$ is the $n$-th power map. We shall identify 
$\mu_n$ with $M_X/nM_X$ according to the identification in [Mochizuki2], 
the discussion at the beginning of \S 2. 

Consider a subset
$$E\subset X^{\cl}.$$
(We will set $E=E_X$ eventually, but $E$ is arbitrary for the present.) 
Let $S\subset X^{\cl}-E$ be a finite set. If we consider the above 
Kummer exact sequence on 
the \'etale site of $U_S \overset \text {def}\to=X-S$ 
and pass to the inverse limit with respect to $n$, then we obtain a  
natural homomorphism
$$\Gamma (U_S,\Cal O_{U_S}^{\times})\to H^1(\Pi _{U_S},M_X)$$
(cf. loc. cit.). (Note that here it suffices to consider the group cohomology of 
$\Pi _{U_S}$ (i.e., as opposed to the \'etale cohomology of $U_S$),
since the extraction of $n$-th roots of an element of $\Gamma (U_S,\Cal O_
{U_S}^{\times})$ yields finite \'etale coverings of $U_S$ that correspond to 
open subgroups of $\Pi _{U_S}$.) 
Observe that this is injective if $\Sigma ^{\dag}=\Primes - \{p\}$, 
since the abelian group $\Gamma (U_S,\Cal O_{U_S}
^{\times})$ is finitely generated and free of $p$-torsion, hence injects into
its pro-prime-to-$p$ completion. 

In particular, by allowing $S$ to vary among
all finite subsets of $X^{\cl}-E$, we obtain a natural 
homomorphism
$$\Cal O_{E}^{\times}\to \underset {S}\to {\varinjlim}\ H^1(\Pi _{U_S},M_X),$$  
where 
$$\Cal O_{E}^{\times}\overset \text {def}\to=\{f\in K_X^{\times}\mid
\sup (\div (f))\cap E=\emptyset\}$$ 
is the multiplicative group of the units
in the 
ring $\Cal O_{E}$ of functions on $X$ 
which are regular at all points of $E$. 
(Here, $K_X$ denotes the function field of $X$.) 
Observe that this is injective
if $\Sigma ^{\dag}=\Primes -\{p\}$.

\proclaim {Proposition 3.1 (Kummer Classes of Functions)} 
Suppose that $S\subset X^{\cl}-E$ is a finite subset. 
Write 
$$\Delta _{U_S}\twoheadrightarrow \Delta _{U_S}^{\text{\rm c-ab}}\twoheadrightarrow 
\Delta _{U_S}^{\text{\rm c-cn}}$$
for the maximal cuspidally abelian and 
the maximal cuspidally central quotients,
respectively, relative to the map $\Delta _{U_S}\twoheadrightarrow \Delta _X$,
and
$$\Pi _{U_S}\twoheadrightarrow \Pi _{U_S}^{\text{\rm c-ab}}\twoheadrightarrow 
\Pi _{U_S}^{\text{\rm c-cn}}$$
for the corresponding quotients of $\Pi _{U_S}$ (i.e., $\Pi_{U_{S}}^{\text{\rm c-ab}}
\defeq \Pi_{U_{S}}/\Ker (\Delta _{U_{S}}\twoheadrightarrow 
\Delta_{U_{S}}^{\text{\rm c-ab}})$,
$\Pi_{U_{S}}^{\text{\rm c-cn}}
\defeq \Pi_{U_{S}}/\Ker (\Delta _{U_{S}}\twoheadrightarrow 
\Delta_{U_{S}}^{\text{\rm c-cn}})$). 
If $x\in X^{\cl}$, then we shall write
$$D_x[U_S]\subset \Pi_{U_S}$$
for the decomposition group at $x$ in $\Pi_{U_S}$ (which is well-defined 
up to conjugation), and $I_x[U_S]  \overset \text {def}\to=D_x[U_S]\cap 
\Delta _{U_S}$ for the inertia subgroup of $D_x[U_S]$. Thus, when 
$x\in S$ we have a natural isomorphism of $M_X$ with 
$I_x[U_S]$ (cf. [Mochizuki2], Proposition 1.6(ii)(iii)). 
Then: 

\noindent
{\rm (i)} The natural surjections induce the following isomorphisms:
$$ H^1(\Pi _{U_S}^{\text{\rm c-cn}},M_X)\isom H^1(\Pi _{U_S}^{\text{\rm c-ab}},M_X)\isom 
H^1(\Pi _{U_S},M_X)$$
In particular, we obtain the following natural homomorphisms: 
$$\Gamma (U_S,\Cal O_{U_S}^{\times})\to H^1(\Pi _{U_S}^{\text{\rm c-cn}},M_X)\isom 
H^1(\Pi _{U_S}^{\text{\rm c-ab}},M_X)\isom
H^1(\Pi _{U_S},M_X),$$
$$\Cal O_{E}^{\times}\to \underset {S}\to {\varinjlim}
\ H^1(\Pi _{U_S}^{\text{\rm c-cn}},M_X)\isom \underset {S}\to 
{\varinjlim}\ H^1(\Pi _{U_S}^{\text{\rm c-ab}},M_X)\isom \underset 
{S}\to {\varinjlim}\ H^1(\Pi _{U_S},M_X),$$
where $S$ varies among all finite subsets of $X-E$. 

These natural homomorphisms are injective, if, moreover, 
$\Sigma^{\dag}=\Primes - \{p\}$.

\noindent
{\rm (ii)} 
Restricting 
cohomology classes of $\Pi _{U_S}$ to the various $I_x[U_S]$ for $x\in S$ 
yields a natural exact sequence:
$$1\to (k^{\times})^{\Sigma^{\dag}}\to H^1(\Pi _{U_S},M_X) \to 
(\underset {s\in S}\to \oplus 
\hat \Bbb Z^{\Sigma^{\dag}})$$
(where we identify $\Hom _{\hat \Bbb Z^{\Sigma^{\dag}}}(I_x[U_S],M_X)$ with $\hat 
\Bbb Z^{\Sigma^{\dag}}$). Moreover, the image (via the natural homomorphism given in
{\rm (i)}) of $\Gamma (U_S,\Cal O_{U_S}^{\times})$ in  $H^1(\Pi _{U_S},M_X)$ is 
equal to the inverse image in  $H^1(\Pi _{U_S},M_X)$ of the submodule of
$$(\underset {s\in S}\to \oplus \Bbb Z)\subset (\underset {s\in S}\to \oplus 
\hat \Bbb Z^{\Sigma^{\dag}})$$
determined by the principal divisors (with support in $S$). A similar 
statement holds when $\Pi _{U_S}$ is replaced by $\Pi _{U_S}^{\text{\rm c-cn}}$ 
or $\Pi _{U_S}^{\text{\rm c-ab}}$.

\noindent
{\rm (iii)} 
If $f\in \Gamma (U_S,\Cal O_{U_S}^{\times})$, write
$$\kappa _f^{\text{\rm c-cn}}\in H^1(\Pi _{U_S}^{\text{\rm c-cn}},M_X),\ \ 
\kappa _f^{\text{\rm c-ab}}\in H^1(\Pi _{U_S}^{\text{\rm c-ab}},M_X),\ \ 
\kappa _f\in H^1(\Pi _{U_S},M_X)$$
for the associated Kummer classes. If $x\in (X^{\cl}-E)-S$, 
then $D_x[U_S]$ maps, via the natural surjection $\Pi _{U_S}\twoheadrightarrow G_k$, 
isomorphically onto the open 
subgroup $G_{k(x)}\subseteq G_k$ (where $k(x)$ is the residue field of $X$ at 
$x$). Moreover, the images of the pulled back classes
$$
\align
\kappa _f^{\text{\rm c-cn}}|_{D_x[U_S]}=\kappa _f^{\text{\rm c-ab}}|_{D_x[U_S]}
=\kappa _f|_{D_x[U_S]}
\in H^1(D_x[U_S],M_X)&\simeq  H^1(G_{k(x)},M_X)\\
&\simeq (k(x)^{\times})^{\Sigma^{\dag}}
\endalign
$$
in $(k(x)^{\times})^{\Sigma^{\dag}}$ are equal to 
the image in $(k(x)^{\times})^{\Sigma^{\dag}}$ of the value 
$f(x)\in k(x)^{\times}$ of $f$ at $x$. 
\endproclaim

\demo {Proof} See [Mochizuki2], Proposition 2.1. (Strictly speaking, 
Proposition 2.1(ii) loc. cit. only treats the case where $S\subset X(k)$, 
but the same proof works well for the general case.) 
\qed
\enddemo

\definition{Remark 3.2} (cf. [Mochizuki2], Remark 2.1.1.) 
In the situation of 
Proposition 3.1(iii), assume $x\in X(k)$ and $S\subset X(k)$ for simplicity.
If we think of the extension $\Pi _{U_S}^{\text{\rm c-cn}}$ 
of $\Pi_X$ as being given by the extension $\Cal D_S$, where $\Cal D$ is a
fundamental extension of $\Pi _{X\times X}$ (cf. Proposition 2.4(i)), 
then it follows that the 
image of $D_x[U_S]$ in $\Pi _{U_S}^{\text{\rm c-cn}}$ may be thought of 
as the image of $D_x[U_S]$ in $\Cal D_S$. 
This image of $D_x[U_S]$ in $\Cal D_S$
amounts to a section of $\Cal D_S\twoheadrightarrow \Pi_X
\twoheadrightarrow G_k$ lying over the section $s_x:G_k\to \Pi_X$ 
corresponding to the rational point $x$ (which is well-defined up to
conjugation). Since $\Cal D_S$ is defined as a 
certain fiber product, this section is equivalent to a collection 
of sections (regarded as ``cyclotomically outer homomorphisms'', i.e., well-defined
up to composition with an inner automorphism of $\Cal D_{y,x}$ by an element
of $\Ker (\Cal D_{y,x}\twoheadrightarrow G_k)$)
$$\gamma _{y,x}:G_k\to \Cal D_{y,x},$$
where $y$ ranges over all points of $S$. Namely, from this point 
of view, Proposition 3.1(iii) may be regarded as saying that the image
in $(k(x)^{\times})^{\Sigma^{\dag}}
=(k^{\times})^{\Sigma^{\dag}}$ of the value $f(x)$ of the function 
$f\in \Gamma (U_S,\Cal O_{U_S}^{\times})$ at $x\in X(k)$ may be 
computed from its Kummer class, as
soon as one knows the sections $\gamma _{y,x}:G_k\to \Cal D_{y,x}$ 
for $y\in S$. Observe that $\gamma _{y,x}$ depends only on $x$, $y$, and 
not on the choice of $S$.
\enddefinition

\definition {Definition 3.3} (cf. [Mochizuki2], Definition 2.2.)
For $x,y\in X(k)$ with $x\neq y$, we shall 
refer to the above section (regarded as a cyclotomically outer homomorphism)
$$\gamma _{y,x}:G_k\to \Cal D_{y,x}$$
as the Green's trivialization of $\Cal D$ at $(y,x)$. If $D$ is a 
divisor on $X$ supported on $k$-rational points $\neq x$, 
then multiplication of the various  Green's trivializations for 
the points in the support of $D$ yields a section (regarded as a 
cyclotomically outer homomorphism)
$$\gamma _{D,x}:G_k\to \Cal D_{D,x}$$
which we shall refer to as the Green's trivialization of $\Cal D$ at 
$(D,x)$.
\enddefinition

\definition {Definition 3.4} (cf. [Mochizuki2], Definition 2.4.) 
Let the notations and the assumptions as in Corollary 2.3. 

\noindent
{\rm (i)} Write $\Cal D$ (respectively, $\Cal E$) for the  
fundamental extension of $\Pi _{X\times X}$ (respectively, 
$\Pi _{Y\times Y}$) that arises as the 
quotient of $\Pi _{U_{X\times X}}^{\text{\rm c-ab}}$ (respectively, $\Pi _{U_
{Y\times Y}}^{\text{\rm c-ab}}$) by the kernel of the maximal cuspidally central
quotient $\Delta _{U_{X\times X}}^{\text{\rm c-ab}}\twoheadrightarrow 
\Delta _{U_{X\times X}}^{\text{\rm c-cn}}$ (respectively,  $\Delta _{U_{Y\times Y}}
^{\text{\rm c-ab}}\twoheadrightarrow \Delta _{U_{Y\times Y}}^{\text{\rm c-cn}}$) 
(cf. Proposition 2.4(ii)). The isomorphism 
$\alpha ^{\text{\rm c-ab}}$ induces naturally an isomorphism:
$$\alpha ^{\text{\rm c-cn}}:\Cal D\isom \Cal E$$
We shall say that $\alpha$ is $(S,T)$-locally Green-compatible 
outside exceptional sets if, for
every pair of points $(x_1,x_2)\in X(k_X)\times X(k_X)$ corresponding via 
$\phi$ to a pair of points $(y_1,y_2)\in Y(k_Y)
\times Y(k_Y)$, such that $x_1\in (X^{\cl}-E_X)-S$, 
$y_1\in (Y^{\cl}-E_Y)-T$, $x_2\in S$, $y_2\in T$, the isomorphism
$$\Cal D_{x_1,x_2}\isom \Cal E_{y_1,y_2}$$
(obtained by restricting $\alpha ^{\text{\rm c-cn}}$ to the various decomposition groups) 
is compatible with the 
Green's trivializations. We shall say that $\alpha$ is 
$(S,T)$-locally principally
Green-compatible outside exceptional sets 
if, for every point $x\in X(k_X)\cap S$ and 
every 
principal divisor
$D$ supported on $k_X$-rational points 
$\neq x$ contained in $X^{\cl}-E_X$, 
corresponding via $\phi$ to a pair $(y,E)$ (so $y\in Y(k_Y)\cap T$), 
the isomorphism
$$\Cal D_{D,x}\isom \Cal E_{E,y}$$
obtained from $\alpha^{\text{\rm c-cn}}$
is compatible with the Green's trivializations.

\noindent
{\rm (ii)} We shall say that $\alpha$ is totally globally 
Green-compatible (respectively, 
totally globally principally Green-compatible) 
outside exceptional sets 
if, for all pair of connected finite \'etale coverings $\xi: X'\to X$, 
$\eta: Y'\to Y$ that arise from open subgroups $\Pi_{X'}\subseteq \Pi_X$, 
$\Pi _{Y'}\subseteq \Pi_Y$, 
corresponding to each other via $\alpha$, then for any subset $S\subset X^{\cl}-E_X$ 
that corresponds, via $\phi$, to $T\subset Y-E_Y$ the isomorphism
$$\Pi_{X'}\isom \Pi_{Y'}$$
induced by $\alpha$ is $(S',T')$-locally Green-compatible (respectively, 
$(S',T')$-locally principally Green-compatible) 
outside exceptional sets, 
where $S'\defeq\xi ^{-1}(S)\subset X^{'\cl}$, 
$T'\defeq\eta^{-1} (T)  \subset Y^{'\cl}$ 
are the inverse images of $S$, $T$, respectively.
\enddefinition

\definition {Remark/Definition 3.5} 
Let $J=J_X$ be the Jacobian variety of $X$. 
Let $\Div _{X- E_X}^0$ be the group of degree zero divisors 
on $X$ which are supported on points in $X-E_X$. 
Write $D_{X- E_X}$ for the kernel of the natural homomorphism 
$\Div _{X- E_X}^0\to J(k)^{\Sigma}$. 
Here, $J(k)^{\Sigma}$ stands for the maximal $\Sigma$-quotient 
$J(k)/(J(k)\{\Sigma'\})$ of $J(k)$, where, 
for an abelian group $M$, $M\{\Sigma'\}$ stands 
for the subgroup of torsion elements $a$ of $M$ such 
that every prime divisor of the order of $a$ 
belongs to $\Sigma'$. 
Then $D_{X- E_X}$ sits naturally in the following exact sequence:
$$0\to \Pri  _{X- E_X}\to D_{X- E_X}\to J(k)\{\Sigma'\} \to 0,$$
where $\Pri _{X- E_X}\defeq \Cal O_{E_X}^{\times}/k^{\times}$ 
stands for the group of principal divisors supported in $X-E_X $.
\enddefinition

\proclaim {Theorem 3.6 (Reconstruction of Functions)} In the situation of 
Theorem 2.2, assume that $\alpha$ is Frobenius-preserving. 
Write $\Sigma\defeq\Sigma_X=\Sigma_Y$ and 
$p\defeq p_X=p_Y$ (cf. Proposition 1.15(ii)(iii)). 
Then:

\noindent
{\rm (i)} The bijection $\phi:X^{\cl}-E_X\isom Y^{\cl}-E_Y$
induced by $\alpha$ (where $E_X$ and $E_Y$ are the exceptional sets) 
induces a natural bijection between the groups 
$D_{X- E_X}$, $D_{Y- E_Y}$. 

\noindent
{\rm (ii)} Assume, moreover, $\Sigma ^{\dag}=\Primes-\{p\}$. 
Then the bijection in (i), together with the isomorphisms 
in Corollary 2.3, induces naturally an injective homomorphism
$$\Cal O_{E_X}^{\times} \hookrightarrow (\Cal O_{E_Y}^{\times})^{p^{-n}},$$
where $p^n$ is the exponent 
of the $p$-primary abelian group $J_Y(k_Y)\{\Sigma'\}$. 

Moreover, this injective homomorphism
$\Cal O_{E_X}^{\times} \hookrightarrow (\Cal O_{E_Y}^{\times})^{p^{-n}}$
is functorial in $X$, $Y$, in the following sense: if $\xi :X'\to X$ is a finite 
\'etale covering, arising from an open subgroup $\Pi_{X'}\subseteq \Pi_X$, which 
corresponds to a finite \'etale covering $\eta:Y'\to Y$ via $\alpha$ (thus, 
$\Pi_{Y'}=\alpha (\Pi_{X'}$)), then we have a commutative diagram:

$$
\CD
\Cal O_{E_{X'}}^{\times} @>>> (\Cal O_{E_{Y'}}^{\times})^{p^{-n'}}\\
@AAA                        @AAA  \\
\Cal O_{E_X}^{\times} @>>>  (\Cal O_{E_Y}^{\times})^{p^{-n}}\\
\endCD
$$
where $E_{X'}\defeq \xi ^{-1}(E_X)$, $E_{Y'}\defeq \eta ^{-1}(E_Y)$, 
$p^{n'}\ge p^n$ is the exponent 
of the $p$-primary abelian group
$J_{Y'}(k_{Y'})\{\Sigma'\}$, and 
the vertical arrows are the natural embeddings.
\endproclaim

\demo {Proof} (cf. [Mochizuki2], Theorem 2.5(ii).)  

\noindent
(i) First, 
$\phi$ induces naturally a bijection  $\Div _{X- E_X}^0
\isom \Div _{Y- E_Y}^0$. 
Second, the natural
homomorphism $\Div _{X- E_X}^0\to J(k_X)^{\Sigma}$ can be recovered
$\varphi$-group-theoretically from $\Pi_X$. 
(Observe the pro-$\Sigma$ version of [Mochizuki2], Proposition 2.3(i). 
See also the discussion before Proposition 2.3 loc. cit.) 
Thus, $\alpha$ induces naturally a commutative diagram:
$$
\CD
\Div _{X- E_X}^0 @>>>  J_X(k_X)^{\Sigma}\\
@VV{\wr}V            @VV{\wr}V  \\
\Div _{Y- E_Y}^0 @>>> J_Y(k_Y)^{\Sigma}\\
\endCD
$$
where the vertical arrows are the isomorphisms induced 
by $\alpha$. From 
this diagram we deduce naturally an isomorphism $D_{X- E_X}
\isom D_{Y- E_Y}$ between the kernels of the horizontal arrows. 

\noindent
(ii) From the isomorphism $D_{X- E_X}\isom D_{Y- E_Y}$ 
we deduce naturally an embedding 
$(\Pri_{X- E_X})^{p^{n}}\hookrightarrow \Pri_{Y- E_Y}$, from which
we deduce an embedding
$(\Cal O_{E_X}^{\times})^{p^{n}}\hookrightarrow (\Cal O_{E_Y}^{\times})$ 
by Corollary 2.3 and Proposition 3.1(i)(ii). 
Finally, the desired commutativity of diagram follows 
easily from the functoriality of Kummer theory.
\qed
\enddemo

\proclaim
{Theorem 3.7 (Totally Globally Principally Green-Compatible Isomorphisms 
Outside Exceptional Sets)} 
In the situation of Theorem 3.6, assume further that 
$\Sigma^{\dag}=\Primes-\{p\}$, and that $\alpha$ is 
totally globally principally Green-compatible outside exceptional sets. 
Then $\alpha$ arises from a uniquely determined commutative diagram of schemes:

$$
\CD
\Tilde X @>{\sim}>> \Tilde Y \\
@VVV   @VVV \\
X @>{\sim}>> Y \\
\endCD
$$
in which the horizontal arrows are isomorphisms and the vertical arrows are the 
profinite \'etale coverings determined by the groups $\Pi_X$, $\Pi_Y$.
\endproclaim

\demo{Proof} 
(cf. [Mochizuki2], Corollary 2.7, Remark 2.7.1.) 
Let $l\neq p$ be a prime number and let $k_X^l$, $k_Y^l$ 
be the (unique) $\Bbb Z_l$-extensions of $k_X$, $k_Y$, respectively. 
Let $X^l$, $Y^l$ be the normalizations of $X$, $Y$ in $K_Xk_X^l$, $K_Yk_Y^l$, 
respectively.
The exponents of the $p$-primary 
abelian groups $J_X(k_X^l)\{\Sigma'\}$, $J_Y(k_Y^l)\{\Sigma'\}$ 
of $J_X(k_X^l)$, $J_Y(k_Y^l)$, respectively, are finite 
(cf. [Rosen], Theorem 11.6). 
So, write $p^{n_0}$ for the exponent of $J_Y(k_Y^l)\{\Sigma'\}$. 
By passing to the limit over the finite extensions of $k_X$, $k_Y$ 
contained in $k_X^l$, $k_Y^l$, respectively, we obtain a natural embedding 
$\Cal O_{E_{X^{l}}}^{\times} \hookrightarrow 
(\Cal O_{E_{Y^{l}}}^{\times})^{p^{-{n_0}}}$, 
where $\Cal O_{E_{X^{l}}}^{\times}$, $\Cal O_{E_{Y^{l}}}^{\times}$ are 
the multiplicative groups of functions on $X^l$, $Y^l$, whose divisor has 
a support disjoint from 
$E_{X^{l}}\defeq E_X\times _{k_X} k_X^l$, 
$E_{Y^{l}}\defeq E_Y\times _{k_Y} k_Y^l$, respectively
(cf. Theorem 3.6(i)). 
(Recall that $E_X\subset X^{\cl}$, $E_Y\subset Y^{\cl}$ are finite 
by Proposition 1.8(vi).) 
The above embedding $\Cal O_{E_{X^{l}}}^{\times} 
\hookrightarrow (\Cal O_{E_{Y^{l}}}^{\times})^{p^{-{n_0}}}$
arises, by Proposition 4.4, from a uniquely determined embedding $K_{X^l}
\hookrightarrow K_{Y^l}^{p^{-{n_0}}}$ of function fields, 
where $K_{X^l}$, $K_{Y^l}$ are the function fields of $X^l$, $Y^l$, respectively. 
(Observe that the value-preserving assumption in Proposition 4.4 
is equivalent to the 
Green-compatibility assumption. See Remark 3.2. 
Observe also that $X^l(k_X^l)$ is an infinite set by the Weil estimate 
on numbers of rational points of curves over finite fields.) 
This embedding of fields 
restricts to the original embedding of multiplicative groups
$\Cal O_{E_{X}}^{\times} \hookrightarrow (\Cal O_{E_{Y}}^{\times})^{p^{-{n_0}}}$ 
(i.e., the restriction of $\Cal O_{E_{X^{l}}}^{\times} \hookrightarrow 
(\Cal O_{E_{Y^{l}}}^{\times})^{p^{-{n_0}}}$), 
hence to an embedding of fields 
$K_X\hookrightarrow K_Y^{p^{-{n_0}}}$. From the 
construction in the proof of Theorem 3.6,  
the image of the emdedding 
$\Cal O_{E_{X}}^{\times} \hookrightarrow (\Cal O_{E_{Y}}^{\times})^{p^{-{n_0}}}$ 
contains $(\Cal O_{E_{Y}}^{\times})^{p^{m_0}}=
((\Cal O_{E_{Y}}^{\times})^{p^{-{n_0}}})^{p^{(n_0+m_0)}}$, where 
$p^{m_0}$ is the exponent of $J_X(k_X^l)\{\Sigma'\}$. Thus, the 
image of the embedding $K_X\hookrightarrow K_Y^{p^{-{n_0}}}$ contains 
$(K_Y^{p^{-{n_0}}})^{p^{n_0+m_0}}=K_Y^{p^{{m_0}}}$. From this, 
we deduce that the embedding
$K_X\hookrightarrow K_Y^{p^{-{n_0}}}$ is radicial and 
maps $K_X$ isomorphically onto $K_Y^{p^{s}}$ 
for some integer $ -n_0\leq s \leq m_0$. 
Thus, in particular, the original embedding 
$\Cal O_{E_{X}}^{\times} \hookrightarrow 
(\Cal O_{E_{Y}}^{\times})^{p^{-{n_0}}}$ 
induces an isomorphism 
$\Cal O_{E_{X}}^{\times} \isom
(\Cal O_{E_{Y}}^{\times})^{p^s}$. Now, it follows from the 
construction in the proof of Theorem 3.6 that 
the image $\text{Im}(\Cal O_{E_{X}}^{\times})$ 
of $\Cal O_{E_{X}}^{\times}$ in 
($\Cal O_{E_{Y^{l}}}^{\times})^{p^{-{n_0}}}$ 
(i.e., $(\Cal O_{E_{Y}}^{\times})^{p^s}$) must be ``commensurate'' 
to $\Cal O_{E_{Y}}^{\times}$ (that is to say, 
the intersection $\text{Im}(\Cal O_{E_{X}}^{\times}) \cap\Cal O_{E_{Y}}^{\times}$
has finite indices both in $\text{Im}(\Cal O_{E_{X}}^{\times})$ 
and in $\Cal O_{E_{Y}}^{\times}$). This implies $s=0$. That is to say, 
the embedding
$K_X\hookrightarrow K_Y^{p^{-{n_0}}}$ is radicial and 
maps $K_X$ isomorphically onto $K_Y$. 

If $\xi :X'\to X$ is a finite \'etale covering, 
arising from an open subgroup $\Pi_{X'}\subseteq \Pi_X$, 
which corresponds to a finite \'etale covering $\eta:Y'\to Y$ via $\alpha$ (thus, 
$\Pi_{Y'}=\alpha (\Pi_{X'}$)), then the commutative diagram in Theorem 3.6(ii), 
together with the above argument,
induces a commutative diagram of embeddings of fields: 

$$
\CD
K_{X'} @>{\tau'}>> K_{Y'}^{p^{-{n_0}'}}\\
@AAA                        @AAA  \\
K_X @>{\tau}>>  K_Y^{p^{-{n_0}}}\\
\endCD
$$
where the vertical arrows are the natural embeddings,  
$p^{{n_0}'}$, $p^{n_0}$ stand for the exponents of the $p$-primary 
abelian groups $J_{Y'}(k_{Y'}^l)$, $J_Y(k_Y^l)$, respectively
(note that $n_0'\ge n_0$),  
and the horizontal arrows are the embeddings obtained above.
Applying the above arguments to $\Pi_X\isom\Pi_Y$ and 
$\Pi_{X'}\isom\Pi_{Y'}$, we obtain $\tau(K_X)=K_Y$, 
$\tau'(K_{X'})=K_{Y'}$. Thus, 
this diagram induces naturally a commutative diagram:


$$
\CD
K_{X'} @>{\sim}>> K_{Y'}\\
@AAA                        @AAA  \\
K_X @>{\sim}>>  K_Y\\
\endCD
$$
where the vertical arrows are the natural embeddings and the horizontal arrows are 
isomorphisms of fields. By passing to the limit over all 
open subgroups of $\Pi_X$ we obtain a 
natural commutative diagram:

$$
\CD
K_{\Tilde X} @>{\sim}>> K_{\Tilde Y}\\
@AAA                        @AAA  \\
K_X @>{\sim}>>  K_Y\\
\endCD
$$
where $K_{\Tilde X}$, $K_{\Tilde Y}$ stand for
the function fields of $\tilde X$, $\Tilde Y$, respectively,  
the vertical arrows are the natural embeddings, and the horizontal arrows 
are isomorphisms of fields. This commutative diagram yields a 
commutative diagram of schemes in the assertion of Theorem 3.7 
with desired properties. (cf. the proof of [Tamagawa1], Theorem (6.3).) 
\qed
\enddemo

\proclaim{Proposition 3.8 (Total Global Green-Compatibility Outside Exceptional Sets)} 
In the situation of Theorem 2.2, assume further that 
$\alpha$ is Frobenius-preserving. 
Then the isomorphism $\alpha$ is totally globally Green-compatible 
outside exceptional sets. In particular, if $\Sigma$ is cofinite, 
then $\alpha$ is totally globally Green-compatible 
outside exceptional sets. 
\endproclaim

\demo
{Proof} For the first assertion, 
the proof of [Mochizuki2], Corollary 3.11 (where 
$\Sigma^{\dag}=\Primes^{\dag}$ and $E_X=E_Y=\emptyset$ are assumed)  
also works well in this case. 
(Thus, the main ingredient of the proof is Theorem 2.7.) 
The second assertion follows from the first, together with 
Proposition 1.15(viii). 
\qed
\enddemo

\proclaim {Theorem 3.9 (A Prime-to-$p$ Version of Grothendieck's 
Anabelian Conjecture
for Proper Hyperbolic Curves over Finite Fields)} Let $X$ and $Y$ be 
proper hyperbolic curves over finite fields $k_X$, $k_Y$, respectively. 
Let $\Sigma_X$, $\Sigma_Y$ be subsets of $\Primes$, and assume 
$\Sigma _X^{\dag}\defeq  \Sigma_X-\{\char (k_X)\}=\Primes-\{\char (k_X)\}$,  
$\Sigma _Y^{\dag}\defeq  \Sigma_Y-\{\char (k_Y)\}=\Primes-\{\char (k_Y)\}$. 
Write $\Pi_X$, $\Pi_Y$ for the geometrically pro-$\Sigma_X$ 
\'etale fundamental group of $X$, the geometrically pro-$\Sigma_Y$ 
\'etale fundamental group of $Y$, respectively. Let
$$\alpha:\Pi_X\isom \Pi_Y$$
be an isomorphism of profinite groups. Then $\alpha$ arises from a uniquely 
determined commutative diagram of schemes:
$$
\CD
\Tilde X @>{\sim}>> \Tilde Y \\
@VVV   @VVV \\
X @>{\sim}>> Y \\
\endCD
$$
in which the horizontal arrows are isomorphisms and the vertical arrows are the 
profinite \'etale coverings corresponding to the groups $\Pi_X$, $\Pi_Y$. 
\endproclaim

\demo
{Proof}\ Follows formally from Theorem 3.7 and Proposition 3.8. 
\qed
\enddemo

As a consequence of Theorem 3.9, we deduce the following:

\proclaim {Corollary 3.10 (A Prime-to-$p$ Version of Grothendieck's Anabelian 
Conjecture for (Not Necessarily Proper) Hyperbolic Curves over Finite Fields)} 
Let $U$, $V$ be (not necessarily proper) hyperbolic curves over 
finite fields $k_U$, $k_V$, 
respectively. 
Let $\Sigma_U$, $\Sigma_V$ be subsets of $\Primes$, and assume 
$\Sigma _U^{\dag}\defeq  \Sigma_U-\{\char (k_U)\}=\Primes-\{\char (k_U)\}$,  
$\Sigma _V^{\dag}\defeq  \Sigma_V-\{\char (k_V)\}=\Primes-\{\char (k_V)\}$. 
Write $\Pi_U$, $\Pi_V$ for the geometrically pro-$\Sigma_U$ 
tame fundamental group of $U$, the geometrically pro-$\Sigma_V$ 
tame fundamental group of $V$, respectively. 
Let
$$\alpha:\Pi_U\isom \Pi_V$$
be an isomorphism of profinite groups. Then $\alpha$ arises from a uniquely 
determined commutative diagram of schemes:
$$
\CD
\Tilde U @>{\sim}>> \Tilde V \\
@VVV   @VVV \\
U @>{\sim}>> V \\
\endCD
$$
in which the horizontal arrows are isomorphisms and the vertical arrows are the 
profinite \'etale coverings corresponding to the groups $\Pi_U$, $\Pi_V$.
\endproclaim

\demo
{Proof} Let $U'\to U$, $V'\to V$ be any 
finite \'etale Galois coverings arising from 
the open normal subgroups $\Pi_{U'}\subseteq \Pi_U$, $\Pi_{V'}\subseteq \Pi_V$, 
which correspond to each other via $\alpha$, such that the smooth compactifications 
$X'$, $Y'$ of $U'$, $V'$, respectively, are hyperbolic, and that the coverings
$\xi:X'\to X$, $\eta : Y'\to Y$, where $X$, $Y$ are the smooth compactifications 
of $U$, $V$, respectively, 
are 
ramified above all the points of $S\defeq X-U$, $T\defeq Y-V$, respectively. 
(Observe that such $U'\to U$, $V'\to V$ are cofinal in 
the finite \'etale coverings arising from open subgroups 
of $\Pi_U$, $\Pi_V$, respectively.)
Thus, the isomorphism $\alpha:\Pi_U\isom\Pi_V$ restricts to 
an isomorphism
$\alpha':\Pi_{U'}\isom \Pi_{V'}$. 
By Proposition 1.15(vi), $\alpha'$ induces naturally an 
isomorphism $\tilde \alpha ':\Pi_{X'}\isom \Pi_{Y'}$, 
which fits into the 
following commutative diagram:
$$
\CD
\Pi_{U'} @>\alpha '>> \Pi_{V'}\\
@VVV            @VVV  \\
\Pi_{X'} @>\tilde \alpha '>> \Pi_{Y'}\\
\endCD
$$
in which the vertical maps are the natural surjections. 

By Theorem 3.9, the isomorphism $\tilde \alpha'$ arises from a 
uniquely determined commutative diagram of schemes:
$$
\CD
\Tilde X' @>{\sim}>> \Tilde Y' \\
@VVV   @VVV \\
X' @>{\sim}>> Y' \\
\endCD
$$
in which the horizontal arrows are isomorphisms and the vertical arrows are the 
profinite \'etale coverings corresponding to the groups $\Pi_{X'}$, $\Pi_{Y'}$.
Since $\tilde\alpha':\Pi_{X'}\isom\Pi_{Y'}$ is 
equivariant with respect to $\alpha:\Pi_U\isom \Pi_V$, 
this last diagram must be also equivariant 
with respet to $\alpha:\Pi_U\isom \Pi_V$. 
Thus, by dividing by the actions of $\Pi_U$, $\Pi_V$, we see that 
it induces naturally a commutative diagram of schemes:
$$
\CD
X' @>{\sim}>> Y' \\
@VV{\xi}V   @VV{\eta}V \\
X @>{\sim}>> Y \\
\endCD
$$
The commutativity of 
this diagram forces the isomorphisms 
$X'\isom Y'$ and $X\isom Y$ to preserve the sets of ramified points 
of $\xi$, $\eta$. Thus, by the choice of $\xi$, $\eta$, this diagram 
induces a commutative diagram of schemes: 
$$
\CD
U' @>{\sim}>> V' \\
@VV{\xi}V   @VV{\eta}V \\
U @>{\sim}>> V \\
\endCD
$$

Finally, by considering this last commutative diagram for any 
coverings $U'\to U$, $V'\to V$ as above, we obtain 
a commutative diagram of schemes in the assertion of Corollary 3.10,  
with desired properties. (cf. the proof of [Tamagawa1], Theorem (6.3).) 
\qed
\enddemo

Finally, we deduce from our main result a prime-to-$p$ birational 
version of Grothendieck's anabelian 
conjecture for function fields of curves over finite fields
(see Corollary 3.11 below).

Let $X$ be a proper, smooth, geometrically connected curve over a finite field $k=k_X$
of characteristic $p=p_X>0$. 
Let $K_X$ be the function field of $X$. Let $G_{K_{X}}\defeq \Gal (K_X^{\sep}/K_X)$ 
be the absolute Galois group of $K_X$ (where $K_X^{\sep}$ is a separable closure of
$K_X$), which sits naturally in the following exact sequence:
$$1\to G_{K_{\overline X}}\to G_{K_{X}} \to G_k\defeq \Gal (\overline k/k)\to 1,$$
where $G_{K_{\overline X}}\defeq \Gal (K_X^{\sep}/K_{\overline X})$  is the absolute 
Galois group of the function field
$K_{\overline X}$ of $\overline X\defeq X\times _{k}\overline {k}$, and 
$G_k\defeq \Gal (\overline k/k)$ is the absolute Galois group of $k$
(here, $\overline k$ is the algebraic closure of $k$ in $K_X^{\sep}$).
Let $G_{K_{\overline X}}'$ be the maximal prime-to-$p$ 
quotient of $G_{K_{\overline X}}$, and let $G_{K_{X}}'\defeq G_{K_X}/\Ker 
(G_{K_{\overline X}}\twoheadrightarrow  G_{K_{\overline X}}')$ 
be the corresponding quotient of 
$G_{K_{X}}$. We shall refer to $G_{K_{X}}'$ as the geometrically 
pro-prime-to-characteristic 
quotient of $G_{K_{X}}$. As an important consequence of Corollary 3.10, we deduce the 
following prime-to-$p$ version of Uchida's Theorem on isomorphisms between absolute 
Galois groups of function fields  (cf. [Uchida]).

\proclaim 
{Corollary 3.11 (A Prime-to-$p$ Version of Uchida's Theorem)} 
Let $X$, $Y$ be proper, smooth, geometrically connected 
curves over finite fields $k_X$, $k_Y$, respectively. 
Let $K_X$, $K_Y$ be the function fields of $X$, $Y$, respectively. Let $G_{K_{X}}$,
$G_{K_{Y}}$ be the absolute Galois groups of $K_X$, $K_Y$, respectively, 
and let $G_{K_{X}}'$, $G_{K_{Y}}'$ be their geometrically 
prime-to-characteristic quotients. Let 
$$\alpha : G_{K_{X}}'\isom G_{K_{Y}}'$$
be an isomorphism of profinite groups. Then $\alpha$ arises from a uniquely 
determined commutative diagram:
$$
\CD
(K_X)\sptilde @>{\sim}>> (K_Y)\sptilde \\
@AAA              @AAA \\
K_{X} @>{\sim}>> K_{Y} \\
\endCD
$$
in which the horizontal arrows are isomorphisms and the vertical arrows are the 
extensions corresponding to the groups $G_{K_{X}}'$,  $G_{K_{Y}}'$, respectively.
\endproclaim

\demo {Proof} 
Following the arguments of [Uchida], Lemma 3 
(involving Brauer groups), 
one can establish a bijection
$\phi:X^{\cl}\isom Y^{\cl}$ such that $\alpha (D_x)=D_{\phi(x)}$ 
holds for each $x\in X^{\cl}$, 
where $D_x$ stands for 
the decomposition group of $G_{K_{X}}'$ at $x$ (which is well-defined
up to conjugation). 
Further, $\alpha (I_x)=I_{\phi (x)}$ also holds for each 
$x\in X^{\cl}$, where $I_x$ stands for the inertia 
subgroup of $D_x$ by 
the same argument (involving local class field theory) 
as in the proof of Lemma 4 loc. cit. 
Let $S\subset X^{\cl}$ be a finite subset such that $U\defeq X-S$ is hyperbolic.
Let $T\defeq \phi (S)$ and $V\defeq Y-T$. Then $\alpha$ induces naturally an isomorphism
$\Pi_U\isom \Pi_V$ between the geometrically prime-to-characteristic quotients of $\pi_1(U)$,
$\pi_1(V)$, respectively. The latter arises, by Corollary 3.10, from a uniquely determined
commutative diagram of schemes:
$$
\CD
\Tilde U @>{\sim}>> \Tilde V \\
@VVV   @VVV \\
U @>{\sim}>> V \\
\endCD
$$
By considering this commutative diagram for all finite subsets 
$S\subset X^{\cl}$, $T\subset Y^{\cl}$ as above, we obtain 
a commutative diagram of field extensions in the assertion of 
Corollary 3.11 with desired properties. 
\qed
\enddemo

\definition {Remark 3.12} In [Stix1], [Stix 2], Stix proved a certain 
relative version 
of Grothendieck's anabelian conjecture 
for hyperbolic curves over finitely generated fields.
His proof relies on (the absolute version of) Grothendieck's anabelian conjecture for 
affine hyperbolic curves over the prime field, proved by Tamagawa in [Tamagawa1].
Using the same arguments as in [Stix1], one should be able to prove a ``prime-to-characteristic''
relative version of Grothendieck's anabelian conjecture for hyperbolic curves over 
finitely generated fields in positive characteristics, by reducing it to our main 
results in Theorem 3.9 and Corollary 3.10. 
\enddefinition

\definition {Remark 3.13}
Even after Theorem 3.9 and Corollary 3.10 are established, it is 
still unclear to the authors, 
at the time of writing, whether or not $E_X=\emptyset$ 
for $\Sigma_X=\Primes-\{\char(k)\}$. 

Indeed, following a standard way in anabelian geometry 
of approaching this kind of problem, let us consider the following tautological 
family of hyperbolic curves of type $(g_X,1)$: 
$$f: U_{X\times X}\defeq X\times X-\iota(X) \to X.$$
Then $f$ induces a right exact sequence: 
$$\Delta_F \to \Delta_{U_{X\times X}} \to \Delta_X \to 1,$$
where $F$ is a geometric fiber of $f$ (which is a hyperbolic 
curve of type $(g_X,1)$), and $\Delta$ stands for the maximal 
pro-$\Sigma_X$ quotient of the geometric fundamental group. 
Suppose that 
this right exact sequence is also left exact: 
$$1 \to \Delta_F \to \Delta_{U_{X\times X}} \to \Delta_X \to 1.$$
Then the sequence 
$$1 \to \Delta_F \to \Pi_{U_{X\times X}} \to \Pi_X \to 1.$$
is also exact, where $\Pi$ stands for the maximal geometrically 
pro-$\Sigma_X$ quotient of the arithmetic fundamental group. 
Now, 
take $x,x'\in X(k)$ and suppose that $D_x, D_{x'}\subset \Pi_X$ 
coincide with each other (up to conjugation). Then, by pulling back the last 
exact sequence by $D_x, D_{x'}\subset \Pi_X$, 
we can easily obtain the following commutative 
diagram: 
$$
\CD
\Pi_{X-\{x\}} @>{\sim}>> \Pi_{X-\{x'\}} \\
@VVV  @VVV \\
\Pi_X @= \Pi_X 
\endCD
$$
Then, by Theorem 3.9 and Corollary 3.10, we obtain 
the following commutative diagram: 
$$
\CD
X-\{x\} @>{\sim}>> X-\{x'\} \\
@VVV  @VVV \\
X @= X. 
\endCD
$$
(Observe that the commutativity follows from the uniqueness 
assertion in Theorem 3.9.) 
This implies $x=x'$, as desired. 

However, 
it is unclear to the authors, 
at the time of writing, 
whether or not the above left exactness (i.e., the injectivity 
of $\Delta_F \to \Delta_{U_{X\times X}}$) is valid. (Note that 
this is a purely topological (or even purely group-theoretical) problem.)
\enddefinition

\subhead
\S 4. Recovering the Additive Structure
\endsubhead

In this \S, we investigate the problem of recovering the additive structure 
of function fields of curves.

Let $X$, $Y$ be proper, smooth, geometrically connected curves over 
fields $k_X$, $k_Y$, respectively. 
Let $X^{\cl}$, $Y^{\cl}$ be the set of closed 
points of $X$, $Y$, respectively. Let $E_X\subset X^{\cl}$, 
$E_Y\subset Y^{\cl}$ be finite subsets, and let
$$\phi:X^{\cl}-E_X\isom Y^{\cl}-E_Y$$
be a (set-theoretic) bijection. 
Write 
$$\Cal O_{E_X}\overset \text {def}\to=\{f\in K_X\mid 
\forall x\in E_X,\  \ord _x(f)\ge 0\},$$
$$\Cal O_{E_Y}\overset \text {def}\to=\{g\in K_Y\mid 
\forall y\in E_Y,\  \ord _y(g)\ge 0\},$$
where $K_X$, $K_Y$ denote the function fields of $X$, $Y$, 
respectively. These are
the semi-local rings of functions on $X$, $Y$ 
that are regular at all points of $E_X$, $E_Y$, respectively. 
Then we have 
$$\Cal O_{E_X}^{\times}=\{f\in K_X^{\times}\mid
\sup (\div(f))\cap E_X=\emptyset\},$$ 
$$\Cal O_{E_Y}^{\times}=\{g\in K_Y^{\times}\mid 
\sup (\div(g))\cap E_Y=\emptyset\}.$$ 
Let
$$\iota : \Cal O_{E_X}^{\times}\hookrightarrow \Cal O_{E_Y}^{\times}$$
be an embedding of multiplicative groups.

\definition
{Definition 4.1} 
The map $\iota : \Cal 
O_{E_X}^{\times}\hookrightarrow \Cal O_{E_Y}^{\times}$ is called 
order-preserving, relative to the bijection $\phi$, if, for 
each $x\in X^{\cl}-E_X$, we have a commutative diagram:
$$
\CD
\Cal O_{E_Y}^{\times}  @>\ord _{\phi (x)}>>   \Bbb Z \\
@A{\iota} AA      @A{e_x}AA \\
\Cal O_{E_X}^{\times}  @>\ord _{x}>>   \Bbb Z \\
\endCD
$$
where the right vertical map is the multiplication by a positive integer
$e_x$ on $\Bbb Z$. 
\enddefinition

\definition {Definition 4.2} 
The map $\iota : \Cal O_{E_X}^{\times}\hookrightarrow \Cal O_{E_Y}
^{\times}$ is called value-preserving, relative to the bijection 
$\phi$, if, for each $x\in X^{\cl}-E_X$, 
there exists an 
embedding of multiplicative groups
$$\iota _x:k(x)^{\times}\hookrightarrow k(\phi (x))^{\times},$$
where $k(x)$, $k(\phi (x))$ are the residue fields of $X$, $Y$ at
$x$, $\phi(x)$, respectively, such that, 
for any $f\in \Cal O_{E_X}^{\times}$ with $\ord_x(f)= 0$, 
$$\ord_{\phi(x)}(\iota(f))=0,\   \iota _x(f(x))=\iota (f)(\phi (x))$$
hold. 
\enddefinition

\definition{Remark 4.3} 
(i) Assume that the map $\iota
$ 
is order-preserving relative to $\phi$. Then $\iota$ induces 
naturally an embedding $k_X^{\times} \hookrightarrow k_Y^{\times}$ 
of the multiplicative groups 
$k_X^{\times}$, $k_Y^{\times}$ of $k_X$, $k_Y$, respectively. 
We extend this embedding to an embedding 
$k_X \hookrightarrow k_Y$ (of multiplicative monoids) by 
$0\mapsto 0$. We denote this last embedding also by $\iota$. 

\noindent
(ii) Assume that the map 
$\iota
$ 
is order-preserving and value-preserving relative to $\phi$. 
For each $x\in X^{\cl}-E_X$, we extend 
$\iota_x
$ 
to an embedding 
$k(x)\cup\{\infty\} \hookrightarrow k(\phi(x))\cup\{\infty\}$ 
by $0\mapsto 0$, $\infty\mapsto\infty$, and denote this last 
embedding also by $\iota_x$. Then the equality 
$$\iota _x(f(x))=\iota (f)(\phi (x))$$
holds for any $x\in X^{\cl}-E_X$ and any $f\in {\Cal O}_{E_X}^{\times}$. 
\enddefinition

Our aim in this \S\  is to prove the following:

\proclaim
{Proposition 4.4 (Recovering the Additive Structure)} Let 
$\iota : \Cal O_{E_X}^{\times}\hookrightarrow \Cal O_{E_Y}^{\times}$
be an embedding of multiplicative groups which is order-preserving and
value-preserving relative to a bijection $\phi :X^{\cl}-
E_X\isom Y^{\cl}- E_Y$, where $E_X\subset X^{\cl}$, $E_Y\subset 
 Y^{\cl}$ are finite subsets. Assume further that 
$X(k_X)$ is an infinite set. Then $\iota$ arises from a uniquely 
determined embedding $K_X \hookrightarrow K_Y$ of function fields.
\endproclaim

The rest of this \S\  will be devoted to the proof of Proposition 4.4. 
Thus, we shall assume that the embedding
$$\iota : \Cal O_{E_X}^{\times}\hookrightarrow \Cal O_{E_Y}^{\times}$$
is order-preserving and value-preserving relative to
a bijection
$$\phi:X^{\cl}-E_X\isom Y^{\cl}-E_Y,$$ 
and that $X(k_X)$ is an infinite set, hence, in particular, 
$k_X$ is an infinite field. 

\proclaim
{Lemma 4.5 (Recovering the Additive Structure of Constants)}
The map $\iota$ preserves the additive structure of the constant 
fields $k_X$, $k_Y$, respectively, i.e., 
$$\iota (\lambda+\mu) = \iota (\lambda)+\iota (\mu)$$ 
holds for any $\lambda, \mu\in k_X$ (cf. Remark 4.3(i)).
\endproclaim

\demo {Proof} 
Fix a point $x_0\in X^{\cl}-E_X$. Then, 
by the Riemann-Roch theorem, we can find a non-constant function
$f\in \Cal O_{E_X}^{\times}$ such that the pole divisor
$\div (f)_{\infty}$ is in the form of $n\cdot x_0$ for 
some integer $n>0$.
Next, observe that $f+\alpha\in \Cal O_{E_X}
^{\times}$ holds for infinitely many $\alpha\in k_X^{\times}$ 
(namely, for any $\alpha \in k_X^{\times}-
(\{-f(x)\mid x\in E_X\}\cap k_X^{\times})$). For 
$\alpha\in k_X^{\times}$ with $f+\alpha\in \Cal O_{E_X}^{\times}$, 
we shall analyze the divisor of the function $\iota (f+\alpha)-\iota (f)$. 
(Observe that $\iota (f+\alpha)-\iota (f)\neq 0$, since 
$\iota$ is injective.) 
We claim: (i) the support of the divisor $\div (\iota (f+\alpha)-\iota (f))$ 
is contained in $\{\phi (x_0)\}\cup E_Y$, and 
(ii) the support of the pole divisor $\div (\iota (f+\alpha)-\iota (f))_\infty$ 
is contained in $\{\phi (x_0)\}$. 
Indeed, let $x\in  X^{\cl}-(E_X\cup \{x_0\})$, and let 
$y=\phi (x)$. Then, $\ord _y(\iota(f+\alpha))=e_x\ord_x(f)\ge 0$ and 
$\ord _y(\iota(f))=e_x\ord_x(f)\ge 0$. 
Moreover, $\iota (f+\alpha)(y)\neq \iota (f)(y)$ as 
follows from the value-preserving assumption, since 
$(f+\alpha)(x)\neq f(x)$. Thus, $y$ does not belong to the support of 
$\div(\iota (f+\alpha)-\iota (f))$,
hence (i) follows. Next, as $\iota(f), \iota(f+\alpha)\in \Cal O_{E_Y}^{\times}$, 
we have $\iota (f+\alpha)-\iota (f)\in \Cal O_{E_Y}$. Thus, 
$\iota (f+\alpha)-\iota (f)$ has no poles in $E_Y$, and (ii) follows. 

Further, the order of $\iota (f+\alpha)-\iota (f)$ 
at the possible pole $\phi (x_0)$ is bounded: 
$$\ord _{\phi (x_0)} (\iota (f+\alpha)-\iota (f))
\geq \min (\ord _{\phi (x_0)} (\iota (f+\alpha)), \ord_{\phi(x_0)}(\iota (f)))
=-n e_{x_0}.$$
We deduce from this 
that there exists an infinite subset $A\subset k_X^{\times}$,
such that $f+\alpha\in \Cal O_{E_X}^{\times}$ holds for all $\alpha\in A$, 
and that the divisor $\div (\iota (f+\alpha)-\iota (f))$ 
is constant for $\alpha\in A$ (i.e., $\div (\iota (f+\alpha)-\iota (f))$ 
($\alpha\in A$) is independent of $\alpha$). 

Let $\alpha, \beta \in A$ with $\alpha\neq\beta$. Thus, 
$$\frac {\iota (f+\beta)-\iota (f)} {\iota (f+\alpha)-\iota (f)}=c \in k_Y^{\times}.$$ 
Further, $c=\frac {\iota (\beta)}{\iota (\alpha)}$, as is easily seen 
by evaluating the function
$\frac {\iota (f+\beta)-\iota (f)} {\iota (f+\alpha)-\iota (f)}$ 
at $\phi (x_1)$, where $x_1$ is a zero of $f$. 
(Observe $x_1\notin E_X$.) 
Thus, we have the equality
$$\iota (\beta) (\iota (f+\alpha)-\iota (f))=\iota (\alpha) (\iota (f+\beta)-
\iota (f))$$ 
which is equivalent to
$$\iota (f)(\iota (\alpha)-\iota (\beta))=\iota (\alpha)\iota (f+\beta)-
\iota (\beta)\iota (f+\alpha).\tag {$*$}$$
Let 
$$g\defeq g_{\alpha,\beta}
\defeq\frac {\beta(f+\alpha)}{(\alpha-\beta)f}\in \Cal O_{E_X}^{\times}.$$ 
Note that $g=\frac{\beta (1+\alpha f^{-1})}{(\alpha-\beta)}$ 
is non-constant, since $f$ is non-constant.
Moreover, we have 
$$g+1=\frac {\alpha(f+\beta)}{(\alpha-\beta)f}\in \Cal O_{E_X}^{\times}.$$ 

We will show the identity $\iota (g+1)=\iota (g)+1$. Indeed, 
$$\iota (g+1)-\iota (g)
=\frac {\iota (\alpha) \iota (f+\beta)}
{\iota (\alpha-\beta) \iota (f)}-\frac {\iota (\beta)\iota (f+\alpha)}
{\iota (\alpha-\beta) \iota (f)}
\overset{(*)}\to=\frac {\iota (\alpha)-\iota (\beta)}
{\iota (\alpha-\beta)}.$$ 
Moreover, 
$$\frac {\iota (\alpha)-\iota (\beta)}{\iota (\alpha-\beta)}=1,$$ 
as follows by evaluating the function $\iota (g+1)-\iota (g)$ at 
$\phi (x_2)$,
where $x_2$ is a zero of $g$. 
(Observe $x_2\notin E_X$.) 
Thus, 
$$\iota (g+1)=\iota (g)+1.$$ 

Next, let $\lambda, \mu\in k_X$, and we shall show the identity 
$\iota (\lambda+\mu)=\iota (\lambda)+\iota (\mu)$. 
If one (or both) of $\lambda$, $\mu$ is $0$, this identity 
clearly holds. So, we may and shall assume 
$\lambda,\mu\in k_X^{\times}$ and 
set $\eta\defeq \lambda/\mu\in k_X^{\times}$. First, 
assume that 
$$\eta \in k_X-(\{g(x)\mid x\in E_X\}\cap k_X)$$
and let $x_3\in X^{\cl}$ be a zero of $g-\eta$. 
Thus, $x_3\notin E_X$, and, by evaluating the identity 
$\iota (g+1)=\iota (g)+1$ at $\phi(x_3)$, we obtain 
$\iota (\eta)+1=\iota (\eta+1)$. 
To show this equality for general $\eta$, 
we shall fix ($f$ and) $\beta\in A$ and make 
$\alpha\in A-\{\beta\}$ vary 
in the expression of $g=g_{\alpha,\beta}$. 
More precisely, take any 
$\alpha\in (A-\{\beta\})-\{\frac{(\eta+1)\beta f(x)}{\eta f(x)-\beta}\mid x\in E_X\}$. 
Then $g=g_{\alpha,\beta}$ satisfies 
$\eta \notin k_X - (\{g(x)\mid x\in E_X\}\cap k_X)$. 
Thus, by the preceding argument, we conclude that 
$$\iota (\eta)+1=\iota (\eta+1)$$ 
holds in general. 

Finally, we obtain
$$\iota (\lambda+\mu)
=\iota (\mu)\iota (\eta +1)
=\iota (\mu) (\iota (\eta)+1)
=\iota (\lambda)+\iota (\mu),$$
as desired. 
\qed
\enddemo

\proclaim{Corollary 4.6}
The map $\iota: k_X\to k_Y$ is an embedding of fields. 
\endproclaim

\demo{Proof}
$\iota$ is multiplicative by definition and additive by Lemma 4.5. 
\qed\enddemo

\proclaim{Corollary 4.7}
For each $x\in X(k_X)-E_X$, the map 
$\iota_x: k(x) \to k(\phi(x))$ is an embedding of fields. 
\endproclaim

\demo{Proof}
For each $x\in X^{\cl}-E_X$, consider the following diagram
$$
\CD
k(x) @>{\iota_x}>> k(\phi(x)) \\
@AAA @AAA \\
k_X @>{\iota}>> k_Y \\
\endCD
$$
where the vertical arrows are evaluation maps. By the value-preserving 
property, this diagram is commutative. If, moreover, $x\in X(k_X)$, 
we have $k_X\isom k(x)$. Thus, Corollary 4.7 follows from Corollary 4.6. 
\qed\enddemo

Next, let $\Bbb Z[\Cal O_{E_X}^{\times}]$, $\Bbb Z[\Cal O_{E_Y}^{\times}]$ 
be the group algebras of the multiplicative groups 
$\Cal O_{E_X}^{\times}$, $\Cal O_{E_Y}^{\times}$, respectively, 
over $\Bbb Z$. The group homomorphism
$\iota:\Cal O_{E_X}^{\times}\hookrightarrow \Cal O_{E_Y}^{\times}$ 
extends uniquely to a ring homomorphism
$$\iota': \Bbb Z[\Cal O_{E_X}^{\times}]\to \Bbb Z[\Cal O_{E_Y}^{\times}].$$
Namely, 
$$\iota' (\sum _i n_if_i)\overset \text {def}\to=\sum _i n_i\iota (f_i)$$
where $n_i\in \Bbb Z$, $f_i\in \Cal O_{E_X}^{\times}$.
Further, let $R_X$, $R_Y$ be the $\Bbb Z$-subalgebras of $K_X$, $K_Y$, 
respectively, generated by $\Cal O_{E_X}^{\times}$, $\Cal O_{E_Y}^{\times}$, 
respectively. Observe that $R_X$, $R_Y$ may be naturally regarded as 
factor rings of $\Bbb Z[\Cal O_{E_X}^{\times}]$, $\Bbb Z[\Cal O_{E_Y}^{\times}]$, 
respectively. 

\proclaim {Lemma 4.8} The ring homomorphism 
$\iota': \Bbb Z[\Cal O_{E_X}^{\times}]\to \Bbb Z[\Cal O_{E_Y}^{\times}]$
induces a (unique) ring homomorphism $\iota_R: R_X\to R_Y$. More precisely, 
The composite of 
$\iota': \Bbb Z[\Cal O_{E_X}^{\times}]\to \Bbb Z[\Cal O_{E_Y}^{\times}]$ 
and the natural surjection $\Bbb Z[\Cal O_{E_Y}^{\times}]\to R_Y$ factors through 
the natural surjection $\Bbb Z[\Cal O_{E_X}^{\times}]\to R_X$. 
\endproclaim

\demo {Proof} 
Take any element 
$$F=\sum _i n_if_i \in \Bbb Z[\Cal O_{E_X}^{\times}],$$
where $n_i\in \Bbb Z$, $f_i\in \Cal O_{E_X}^{\times}$, such that the image 
$F_X$ of $F$ in $R_X$ is $0$. Then we have to show that 
the image $F_Y$ of $F$ in $R_Y$ is also $0$. To avoid confusion, 
we shall denote a sum in a ring $R$ by means of 
$\sum{}_R$. Then the assumption 
$F_X=0$ can be expressed as the equality 
$$\sum_i{}_{R_X}\  n_if_i =0.$$
Let $S_i\subset X^{\cl}$ denote the 
(finite) set of poles of $f_i$ and consider a 
point $x\in X(k_X)-(E_X\cup\cup_i S_i)$. 
By evaluating the above equality at $x$, we obtain 
the equality 
$$\sum _i{}_{k(x)}\  n_if_i(x) =0. $$
By Corollary 4.7 and the value-preserving property at $x$, 
this last equality implies the equality
$$\sum _i{}_{k(\phi(x))}\  n_i\iota(f_i)(\phi(x)) =0,$$
or, equivalently, the equality 
$$F_Y(\phi(x))=0$$
in $k(\phi(x))$. Since $x$ is an arbitrary point in 
the infinite set $X(k_X)-(E_X\cup\cup_i S_i)$, this implies 
$F_Y=0$ in $R_Y$, as desired. 
\qed
\enddemo

\proclaim {Lemma 4.9} We have $R_X=\Cal O_{E_X}$ and
$R_Y=\Cal O_{E_Y}$.
\endproclaim

\demo {Proof} It suffices to prove $R_X=\Cal O_{E_X}$. 
It is clear from the definitions that $R_X\subset \Cal O_{E_X}$.
We will show the other inclusion $\Cal O_{E_X}\subset R_X$. 
Let $f \in\Cal O_{E_X}$, then $(f(x))_{x\in E_X}\in 
\oplus _ {x\in E_X}k(x)$. Fix $\alpha\in k_X^{\times}- \{1\}$. 
For each $x\in E_X$, set 
$$\epsilon (x)=
\cases
f(x)-1, &\text {if $f(x)\neq 1$},\\
1-\alpha, &\text {if $f(x)= 1$}.
\endcases
$$
Thus, we have 
$$(\epsilon (x))_{x\in E_X}\in \oplus _ {x\in E_X}k(x)^{\times}.$$ 
Further, set
$$(\delta (x))_{x\in E_X}\overset \text {def}\to=
(f(x))_{x\in E_X}-(\epsilon (x))_{x\in E_X},$$
then we have 
$$\delta(x)=
\cases
1, &\text {if $f(x)\neq 1$},\\
\alpha, &\text {if $f(x)= 1$}.
\endcases
$$
Thus, in particular, we have 
$$(\delta (x))_{x\in E_X}\in \oplus _ {x\in E_X}k(x)^{\times}.$$ 
Since the evaluation map $\Cal O_{E_X}\to\oplus _ {x\in E_X}k(x)$ is 
surjective, we can find a function $g\in \Cal O_{E_X}$, 
such that $(g(x))_{x\in E_X}=(\epsilon(x))_{x\in E_X}$. 
Set $h\defeq f-g$. Then we have 
$(h(x))_{x\in E_X}=(\delta(x))_{x\in E_X}$. 
Since $g,h\in \Cal O_{E_X}^{\times}$, 
the equality $f=g+h$ shows that $f\in R_X$.
\qed
\enddemo

\proclaim {Lemma 4.10} The ring homomorphism $\iota_R: R_X\to R_Y$ 
in Lemma 4.8 is injective.
\endproclaim

\demo
{Proof} Take any $f\in R_X$ with $\iota_R (f)=0$. 
As in the proof of Lemma 4.9, $f\in R_X=\Cal O_{E_X}$ can be 
written as 
$f=g+h$ for some $g,h\in\Cal O_{E_X}^{\times}$. Now we have 
$$\iota(g)=\iota_R(g)=\iota_R(f)+\iota_R(-h)=\iota(-h).$$
Since $\iota$ is injective, this shows $g=-h$, 
hence $f=0$, as desired. 
\qed
\enddemo

\proclaim
{Corollary 4.11} The ring homomorphism $\iota_R :R_X\to R_Y$ 
extends uniquely to an embedding $K_X\hookrightarrow K_Y$ of fields.
\endproclaim

\demo{Proof} 
This follows from Lemmas 4.9 and 4.10. 
\qed\enddemo

This completes the proof of Proposition 4.4.
\qed

\definition{Remark 4.12}
The above proof of Proposition 4.4 relies 
on the value-preserving property at all but 
finitely many points of $X^{\cl}$ 
(or, more precisely, all points of 
$X^{\cl}-E_X$). This is a contrast to 
the proof of [Tamagawa1], Lemma (4.7), 
which relies on the value-preserving property 
at only finitely many points. Thus, unlike 
the case of [Tamagawa1], Theorem (4.3), 
we need, at least for the time being, 
to resort to Mochizuki's theory of cuspidalizations 
to prove Corollary 3.10, even in the affine case. 
\enddefinition

$$\text{References.}$$

\noindent
[Boxall] Boxall, J., 
Autour d'un probl\`eme de Coleman, 
C. R. Acad. Sci. Paris 315 (1992), S\'er. I, 1063--1066.

\noindent
[Chevalley] Chevalley, C., 
Deux th\'eor\`emes d'arithm\'etique, J. Math. Soc. Japan, 3 (1951), 
36--44. 

\noindent
[Grothendieck] Grothendieck, A., Brief an G. Faltings. (German) [Letter to G. Faltings] 
With an English translation on pp. 285--293. London Math. Soc. Lecture Note Ser., 242,  
Geometric Galois actions, 1,  49--58, Cambridge Univ. Press, Cambridge, 1997.

\noindent
[Mochizuki1] Mochizuki, S., 
The local pro-$p$ anabelian geometry of curves, 
Inventiones Mathematicae 138 (1999), 
319--423. 

\noindent
[Mochizuki2] ---, Absolute anabelian cuspidalizations of proper 
hyperbolic curves, preprint RIMS-1490 (2005), revised version (2007), 
available at 
http://www.kurims.
\linebreak
kyoto-u.ac.jp/\~{}motizuki/papers-english.html, 
to appear in Journal of Mathematics of Kyoto University. 

\noindent
[Nakamura1] Nakamura, H., 
Galois rigidity of the \'etale fundamental groups of punctured projective lines, 
J. Reine Angew. Math., 411 (1990), 205--216. 

\noindent
[Nakamura2] ---, Galois rigidity of algebraic mappings 
into some hyperbolic varieties, 
International Journal of Mathematics, 4 (1993), 421--438. 

\noindent
[Nakamura3] ---, 
Galois rigidity of pure sphere braid groups and profinite calculus, 
J. Math. Sci. Univ. Tokyo, 1 (1994), 
71--136. 

\noindent
[NT] Nakamura, H. and Tsunogai, H., Some finiteness theorems 
on Galois centralizers in pro-$l$ mapping class groups, 
J. Reine Angew. Math., 441 (1993), 115--144.

\noindent
[Rosen] Rosen, M., Number Theory in Function Fields, 
Graduate Texts in Mathematics, 210. Springer-Verlag, New York, 2002.

\noindent
[Stix1] Stix, J., Affine anabelian curves in positive characteristic.  
Compositio Math.  134  (2002),  no. 1, 75--85.

\noindent
[Stix2] ---, Projective anabelian curves in positive characteristic 
and descent theory for log-\'etale covers. Dissertation, Rheinische 
Friedrich-Wilhelms-Universit\"at Bonn, Bonn, 2002. Bonner Mathematische Schriften, 
354. Universit\"at Bonn, 
Mathematisches Institut, Bonn, 2002.

\noindent
[Tamagawa1] Tamagawa, A., 
The Grothendieck conjecture for affine curves, 
Compositio Math. 109 (1997), 
135--194. 

\noindent
[Tamagawa2] ---, On the tame fundamental groups of curves over algebraically closed 
fields of characteristic $>0$, in Galois groups and fundamental groups, 
Math. Sci. Res. Inst. Publ. 41, Cambridge Univ. Press, Cambridge, 
2003, 47--105. 

\noindent
[Tamagawa3] ---, Finiteness of isomorphism classes of curves in positive 
characteristic 
with prescribed fundamental groups, 
J. Algebraic Geom. 13 (2004), 675--724. 

\noindent
[Uchida] Uchida, K., Isomorphisms of Galois groups of algebraic function fields.  
Ann. Math. (2)  106  (1977), no. 3, 589--598. 

\bigskip
\noindent
Mohamed Sa\"\i di

\noindent
School of Engineering, Computer Science, and Mathematics

\noindent
University of Exeter

\noindent
Harrison Building

\noindent
North Park Road

\noindent
EXETER EX4 4QF 

\noindent
United Kingdom

\noindent
M.Saidi\@exeter.ac.uk

\bigskip
\noindent
Akio Tamagawa

\noindent
Research Institute for Mathematical Sciences

\noindent
Kyoto University

\noindent
KYOTO 606-8502

\noindent
Japan

\noindent
tamagawa\@kurims.kyoto-u.ac.jp
\enddocument